# STABLE HOMOLOGY OF SURFACE DIFFEOMORPHISM GROUPS MADE DISCRETE

SAM NARIMAN

ABSTRACT. We answer affirmatively a question posed by Morita [Mor06] on homological stability of surface diffeomorphism groups made discrete. In particular, we prove that $C^\infty$-diffeomorphisms of surfaces as family of discrete groups exhibit *homological stability*. We show that the stable homology of $C^\infty$-diffeomorphims of surfaces as discrete groups is the same as homology of certain infinite loop space related to Haefliger's classifying space of foliations of codimension 2. We use this infinite loop space to obtain new results about (non)triviality of characteristic classes of flat surface bundles and codimension 2 foliations.

## 0. STATEMENTS OF THE MAIN RESULTS

This paper is a continuation of the project initiated in [Nar14] on the homological stability and the stable homology of discrete surface diffeomorphisms.

### 0.1. Homological stability for surface diffeomorphisms made discrete.

To fix some notations, let $\Sigma_{g,n}$ denote a surface of genus $g$ with $n$ boundary components and let $\mathrm{Diff}^\delta(\Sigma_{g,n}, \partial)$ denote the discrete group of orientation preserving diffeomorphisms of $\Sigma_{g,n}$ that are supported away from the boundary.

The starting point of this paper is a question posed by Morita [Mor06, Problem 12.2] about an analogue of Harer stability for surface diffeomorphisms made discrete. In light of the fact that all known cohomology classes of $\mathrm{BDiff}^\delta(\Sigma_g)$ are stable with respect to $g$, Morita [Mor06] asked

**Question** (Morita). *Do the homology groups of $\mathrm{BDiff}^\delta(\Sigma_g)$ stabilize with respect to $g$?*

In order to prove homological stability for a family of groups, it is more convenient to have a map between them. To define a map inducing homological stability, let $j : \Sigma_{g,1} \hookrightarrow \Sigma_{g+1,1} \backslash \partial \Sigma_{g+1,1}$ be an embedding such that the complement of $j(\Sigma_{g,1})$ in $\Sigma_{g+1,1} \backslash \partial \Sigma_{g+1,1}$ is diffeomorphic to the interior of $\Sigma_{1,2}$. By extending diffeomorphisms via the identity, this embedding induces a group homomorphism between diffeomorphism groups $s : \mathrm{Diff}^\delta(\Sigma_{g,1}, \partial) \to \mathrm{Diff}^\delta(\Sigma_{g+1,1}, \partial)$. Although the stabilization map $s$ depends on the embedding $j$, it is not hard to see that different choices of embeddings induce the same map on the group homology (see [Nar14, Theorem 2.5] for more details); hence by abuse of notation, we denote the induced map between classifying spaces $s : \mathrm{BDiff}^\delta(\Sigma_{g,1}, \partial) \to \mathrm{BDiff}^\delta(\Sigma_{g+1,1}, \partial)$ by the same letter. Our first theorem affirmatively answers Morita's question .

**Theorem 0.1.** *The stabilization map*

$$s_* : H_k(\mathrm{BDiff}^\delta(\Sigma_{g,1}, \partial); \mathbb{Z}) \to H_k(\mathrm{BDiff}^\delta(\Sigma_{g+1,1}, \partial); \mathbb{Z})$$

*induces isomorphism as long as $k \leq (2g-2)/3$.*

*Remark* 0.2. Bowden [Bow12] proved stability for $k \leq 3$ as $g \geq 8$. Here, we give a proof with the same stability range as that of the mapping class groups.

These homological stability results hold for surface diffeomorphisms with any order of regularity i.e. the stabilization map induces a homology isomorphism for $C^r$ diffeomorphisms of surfaces as $r > 0$. However, the remarkable theorem of Tsuboi [Tsu89] implies that the classifying space of $C^1$-diffeomorphisms with discrete topology, $\mathrm{BDiff}^{\delta,1}(\Sigma_{g,n}, \partial)$, is homology equivalent to the classifying space of the mapping class group of $\Sigma_{g,n}$. Hence, for regularity $r = 1$, the homological stability is already implied by Harer stability [Har85] for mapping class groups.

It should be further noted that the proof of [Nar14, Theorem 1.1] which is a similar theorem for high dimensional analogue of surfaces, does not carry over verbatim to prove homological stability of surface diffeomorphisms. In proving homological stability for a family of groups, one key step is to build a highly connected simplicial complex on which the family acts. To prove the highly connectedness of the simplicial complex used in [Nar14, Theorem 1.1], it is essential to work in dimension higher than 5 so that certain surgery arguments work.

However, Randal-Williams in [RW16] proved a homological stability theorem for moduli spaces of surfaces equipped with a "tangential structure". We use Thurston's generalization of Mather's theorem [Thu74] in foliation theory and Randal-Williams' theorem ([RW16, Theorem 7.1]) to establish homological stability of $\mathrm{Diff}^\delta(\Sigma_{g,1}, \partial)$. The advantage of this high-powered approach is that it describes the limiting homology in terms of an infinite loop space related to codimension 2 foliations.

### 0.2. Stable homology of $\mathrm{Diff}^\delta(\Sigma_{g,n}, \partial)$.

Analogous to [Nar14, Theorem 1.2], we can describe the stable homology of $\mathrm{BDiff}^\delta(\Sigma_{g,n}, \partial)$ in terms of an infinite loop space related to the Haefliger category. Let us recall the definition of the Haefliger classifying space of foliations.

**Definition 0.3.** The Haefliger category $\Gamma_n^r$ is a topological category whose objects are points in $\mathbb{R}^n$ with its usual topology and morphisms between two points, say $x$ and $y$, are germs of $C^r$ diffeomorphisms that send $x$ to $y$. The space of morphisms is equipped with the sheaf topology (see Section 1.2.1 for more details). If we do not decorate the Haefliger category with $r$, we usually mean the space of morphisms has the $C^\infty$-regularity. By $\mathrm{S}\Gamma_n^r$, we mean the subcategory of $\Gamma_n^r$ with the same objects, but the morphisms are germs of orientation preserving diffeomorphisms (See [Hae71] for more details).

The classifying space of the Haefliger category classifies Haefliger structures up to concordance. The normal bundle to the Haefliger structure induces a map

$$\nu : \mathrm{B}\mathrm{S}\Gamma_n \longrightarrow \mathrm{BGL}_n^+(\mathbb{R}),$$

where $\mathrm{GL}_n^+(\mathbb{R})$ is the group of real matrices with positive determinants.

Let $\gamma_2$ denote the tautological bundle over $\mathrm{GL}_2^+(\mathbb{R})$. Recall that the *Madsen-Tillmann spectrum* $\mathrm{MTSO}(2)$ is the Thom spectrum of the virtual bundle $-\gamma_2$ over $\mathrm{BGL}_2^+(\mathbb{R})$. Let $\mathrm{MT}\nu$ denote the Thom spectrum of the virtual bundle $\nu^*(-\gamma_2)$ over $\mathrm{B}\mathrm{S}\Gamma_2$ (see Definition 2.2 for a more detailed description). We denote the base point component of the infinite loop space associated to this spectrum by $\Omega_0^\infty \mathrm{MT}\nu$. As we shall explain in Section 2.2, there exists a *parametrized Pontryagin-Thom construction* which induces a continuous map

$$\alpha : \mathrm{BDiff}^\delta(\Sigma_{g,n}, \partial) \longrightarrow \Omega_0^\infty \mathrm{MT}\nu,$$

Our second theorem is an analogue of the Madsen-Weiss theorem [MW07] for discrete surface diffeomorphisms.

**Theorem 0.4.** *The map $\alpha$ induces a homology isomorphism in the stable range of Theorem 0.1.*

For any topological group $G$, let $G^\delta$ denote the same group with the discrete topology. The identity map defines a continuous homomorphism $G^\delta \to G$. Thus, for the topological group $\mathrm{Diff}(\Sigma_{g,n}, \partial)$ with $C^\infty$ topology, the identity induces a map

$$\iota \colon \mathrm{BDiff}^\delta(\Sigma_{g,n}, \partial) \to \mathrm{BDiff}(\Sigma_{g,n}, \partial).$$

To study the effect of this map on cohomology in the stable range, we study the natural map that is induced by $\nu$ between the following two infinite loop spaces

$$\Omega^\infty \boldsymbol{\nu} \colon \Omega_0^\infty \mathrm{MT}\nu \to \Omega_0^\infty \mathrm{MTSO}(2).$$

*Remark* 0.5. In foliation theory, Mather and Thurston studied the homotopy fiber of $\iota$ (see Theorem 1.3), and they already translated understanding the homotopy fiber of $\iota$ to a homotopy theoretic question about $\mathrm{B}\overline{\mathrm{S}\Gamma}_2$. However, we pursue a different approach by replacing $\mathrm{BDiff}^\delta(\Sigma_{g,n}, \partial)$ and $\mathrm{BDiff}(\Sigma_{g,n}, \partial)$ with appropriate infinite loop spaces to study $\iota$ in the stable range. As we shall see in section 0.3 and section 0.4, this approach is much more amenable to actual calculations in the stable range instead of understanding the homotopy fiber of $\iota$.

**Theorem 0.6.** *The map $\Omega^\infty \boldsymbol{\nu}$ induces an injection on $\mathbb{F}_p$-cohomology i.e.*

$$H^k(\Omega_0^\infty \mathrm{MTSO}(2); \mathbb{F}_p) \hookrightarrow H^k(\Omega_0^\infty \mathrm{MT}\nu; \mathbb{F}_p),$$

*for any $k$.*

**Corollary 0.7.** *The map $\iota$ induces an injection*

$$H^k(\mathrm{BDiff}(\Sigma_{g,n}, \partial); \mathbb{F}_p) \hookrightarrow H^k(\mathrm{BDiff}^\delta(\Sigma_{g,n}, \partial); \mathbb{F}_p),$$

*as long as $k \leq (2g-2)/3$.*

*Remark* 0.8. All three theorems 0.1, 0.4 and 0.6 in fact hold for $C^r$-diffeomorphisms for any $r > 0$. For applications in Section 0.3 and Section 0.4, we formulated the theorems for $C^\infty$-diffeomorphisms but in fact all of them hold for $C^r$-diffeomorphisms while $r > 1$ and $r \neq 3$. As we mentioned earlier, the case of $\mathrm{Diff}^{\delta,1}(\Sigma_{g,n}, \partial)$ is an exception that thanks to Tsuboi's theorem ([Tsu89]) this group has the same homology of the mapping class group of $\Sigma_{g,n}$. The reason we also exclude $r = 3$ is that for most of the applications, we need the perfectness of the identity component $\mathrm{Diff}_0^{\delta,r}(\Sigma_{g,n}, \partial)$ which in the smooth case is a consequence of Thurston's work ([Thu74]) and for $r \neq 3$ is a consequence of Mather's work ([Mat74]).

0.3. **Applications to characteristic classes of flat surface bundles.** The theory of characteristic classes of fiber bundles and foliated fiber bundles (i.e. fiber bundle with a foliation transverse to fibers) whose fibers are diffeomorphic to a $C^\infty$-manifold M is equivalent to understanding the cohomology groups $H^*(\mathrm{BDiff}(\mathrm{M}))$ and $H^*(\mathrm{BDiff}^\delta(\mathrm{M}))$ respectively. The theory of characteristic classes of manifold bundles and surface bundles in particular have been studied extensively (see [GRW14], [Mor87] and references thereof). Therefore, we have some understanding of $H^*(\mathrm{BDiff}(\mathrm{M}))$ for certain classes of manifolds. For foliated (flat) manifold bundles, however, there seems to be very little known about the existence of non trivial characteristic classes in $H^*(\mathrm{BDiff}^\delta(\mathrm{M}))$.

The abstract results in Section 0.2 shed new light on (non)triviality of characteristic classes of flat surface bundles. They also provide a unified approach to previous results of Morita, Kotschick ([KM07]) and Bowden ([Bow12]). We pursue the study of $H^*(\mathrm{BSymp}^\delta(\Sigma_g))$ in light of these theorems elsewhere.

Morita in [Mor87] showed that finite index subgroups of mapping class group of surfaces cannot be realized as subgroups of diffeomorphisms by showing that Mumford-Miller-Morita classes $\kappa_i \in H^{2i}(\mathrm{BDiff}(\Sigma_g); \mathbb{Q})$ for $i > 2$ get sent to zero via the induced map

$$H^*(\mathrm{BDiff}(\Sigma_g); \mathbb{Q}) \to H^*(\mathrm{BDiff}^\delta(\Sigma_g); \mathbb{Q}).$$

Unlike the cohomology with rational coefficients, Corollary 0.7 implies that all monomials of $\kappa_i$'s are non-torsion classes in $H^*(\mathrm{BDiff}^\delta(\Sigma_g); \mathbb{Z})$.

**Theorem 0.9.** *There is an injection*

$$\mathbb{Z}[\kappa_1, \kappa_2, \dots] \hookrightarrow H^k(\mathrm{BDiff}^\delta(\Sigma_g); \mathbb{Z}),$$

*as long as* $k \leq (2g - 2)/3$.

*Remark* 0.10. Note that not all non-torsion classes in $H^*(\mathrm{BDiff}^\delta(\Sigma_g); \mathbb{Z})$ can be realized by an element in $\mathrm{Hom}(H_*(\mathrm{BDiff}^\delta(\Sigma_g); \mathbb{Z}); \mathbb{Z})$. Because $\mathrm{BDiff}^\delta(\Sigma_g)$ is not a finite type space, the universal coefficient theorem implies that the Ext term in $H^*(\mathrm{BDiff}^\delta(\Sigma_g); \mathbb{Z})$ might have non-torsion classes too. In particular, nontriviality of $\kappa_i$ in $H^*(\mathrm{BDiff}^\delta(\Sigma_g); \mathbb{Z})$ does not imply that there exists a flat surface bundle whose $\kappa_i$ is nonzero. But in fact one can use the method of [AKU01] to prove such flat surface bundles exist.

**Corollary 0.11.** *The group* $H_{2k-1}(\mathrm{BDiff}^\delta(\Sigma_g); \mathbb{Z})$ *is not finitely generated as long as* $k > 2$ *and* $k \leq 2g/3$.

Morita and Kotschick [KM05] constructed a flat surface bundle over a surface whose signature is nonzero. Hence, they conclude that $\kappa_1$ that is 3 times the signature of the total space, is nonzero in $H^2(\mathrm{BDiff}^\delta(\Sigma_g); \mathbb{Q})$. We can use Theorem 0.4 to give a homotopy-theoretic proof of their result in the stable range.

**Theorem 0.12** (Morita-Kotschick). *The image of* $\kappa_1^n$ *in* $H^{2n}(\mathrm{BDiff}^\delta(\Sigma_g); \mathbb{Q})$ *is nonzero for all positive integer* $n$ *provided that* $g \geq 3n$.

To summarize the (non)vanishing results of MMM-classes for flat surface bundles, recall that the Bott vanishing theorem implies that $\kappa_i$ for $i > 2$ vanishes in $H^{2i}(\mathrm{BDiff}^\delta(\Sigma_g); \mathbb{Q})$ and Theorem 0.12 implies that the first MMM-class $\kappa_1$ does not vanish in $H^2(\mathrm{BDiff}^\delta(\Sigma_g); \mathbb{Q})$. Moreover, Theorem 0.9 implies that all $\kappa_i$'s are non-torsion classes in $H^{2i}(\mathrm{BDiff}^\delta(\Sigma_g); \mathbb{Z})$. By Theorem 0.9 we know that $\kappa_2$ is nonzero in $H^4(\mathrm{BDiff}^\delta(\Sigma_g); \mathbb{Z})$, however, we still do not know the answer to the following problem posed by Morita and Kotschick [KM05]:

**Problem.** *Determine whether the second MMM class* $\kappa_2$ *is non-trivial in* $H^4(\mathrm{BDiff}^\delta(\Sigma_g); \mathbb{Q})$.

We prove that this problem is equivalent to an open problem in foliation theory related to the cube of the Euler class of the normal bundle of codimension 2 foliations [Hur03, Problem 15.4].

**Theorem 0.13.** *The MMM-class* $\kappa_2$ *in* $H^4(\mathrm{BDiff}^\delta(\Sigma_g); \mathbb{Q})$ *is nonzero for* $g > 6$ *if and only if a* $C^2$*-foliation* $\mathcal{F}$ *of codimension* 2 *on a* 6 *manifold exists such that* $e(\nu(\mathcal{F}))^3 \neq 0$ *where* $\nu(\mathcal{F})$ *is the normal bundle of the foliation* $\mathcal{F}$ *and* $e(\nu(\mathcal{F}))$ *is its Euler class.*

*Remark* 0.14. Using the universal coefficient theorem, Thom's result on representing cycles by manifolds and Thurston's h-principle for foliations of codimension greater than 2, one can show that the existence of a codimension 2 foliation $\mathcal{F}$ with $e(\nu(\mathcal{F}))^3 \neq 0$ is in fact equivalent to proving $\nu^*(e^3) \in H^6(\mathrm{B}\Gamma_2^+; \mathbb{Q})$ does not vanish where $e \in H^2(\mathrm{BGL}_2^+(\mathbb{R}); \mathbb{Q})$ is the universal Euler class. We then show that nonvanishing of $\kappa_2$ and $\nu^*(e^3)$ are equivalent.

Theorem 0.4 and well-known results about the continuous variation of foliations of codimension 2 can be used to construct more nontrivial classes on flat surface bundles. For example, Rasmussen [Ras80] showed that the two Godbillon-Vey classes $h_1.c_2$ and $h_1.c_1^2$ in $H^5(B\Sigma\Gamma_2;\mathbb{R})$ (see [Bot72, Section 10] for definitions of Godbillon-Vey classes) continuously vary for families of foliations of codimension 2, i.e. the map

$$H_5(B\Sigma\Gamma_2;\mathbb{Q}) \xrightarrow{(h_1c_1^2, h_1c_2)} \mathbb{R}^2,$$

induced by the evaluation of $h_1.c_2$ and $h_1.c_1^2$ is surjective. We use this theorem of Rasmussen to simplify the proof of Bowden's theorem ([Bow12]) that says for all $g$, the fiber integration of the two Godbillon-Vey classes $h_1.c_2$ and $h_1.c_1^2$ induce a surjective homomorphism

$$(0.15) \qquad H_3(\mathrm{BDiff}^\delta(\Sigma_g);\mathbb{Q}) \longrightarrow \mathbb{R}^2.$$

In the stable range though, we use Theorem 0.4 to prove a stronger result that essentially all secondary classes in $H^3(\mathrm{BDiff}^\delta(\Sigma_g);\mathbb{R})$ come from secondary classes of flat disk bundles, more precisely

**Theorem 0.16.** *Let $\mathbb{R}^2 \hookrightarrow \Sigma_g$ be an embedding of an open disk in the surface $\Sigma_g$. For $g \geq 6$, the induced map*

$$H_3(\mathrm{BDiff}_c^\delta(\mathbb{R}^2);\mathbb{Q}) \longrightarrow H_3(\mathrm{BDiff}^\delta(\Sigma_g);\mathbb{Q}),$$

*is surjective.*

Using the surjectivity of the map 0.15 and the Hopf algebra structure on $H_*(\Omega_0^\infty \mathrm{MT}\nu;\mathbb{Q})$, we construct *discontinuous* classes (see [Mor85] for applications of discontinuous invariants) in $H_{3k}(\mathrm{BDiff}^\delta(\Sigma_g);\mathbb{Q})$.

**Theorem 0.17.** *There exists a surjective map*

$$H_{3k}(\mathrm{BDiff}^\delta(\Sigma_g);\mathbb{Q}) \longrightarrow \wedge_{\mathbb{Q}}^k \mathbb{R}^2,$$

*provided $k \leq (2g-2)/9$ where $\wedge_{\mathbb{Q}}^k \mathbb{R}^2$ is the $k$-th exterior power of $\mathbb{R}^2$ as the vector space over $\mathbb{Q}$ .*

0.4. **Applications to the foliated cobordism of codimension 2.** Let $\mathcal{F}\Omega_{n,k}$ be the cobordism group of $n$-manifolds with a foliation of codimension $k$ and $\mathrm{MSO}_n(X)$ be the oriented cobordism group of $n$-manifolds equipped with a map to $X$. Using Theorem 0.4, we compare codimension 2 foliations with foliated surface bundles. By Atiyah-Hirzebruch spectral sequence, the homology equivalence in Theorem 0.4 implies the following equivalence in bordism theory

$$\mathrm{MSO}_n(\mathrm{BDiff}^\delta(\Sigma_g)) \longrightarrow \mathrm{MSO}_n(\Omega_0^\infty \mathrm{MT}\nu),$$

which is an isomorphism in the stable range. Let $e_n$ be that map that associates to every flat surface bundle, the foliated cobordism class of the codimension 2 foliation on the total space of the surface bundle

$$e_n : \mathrm{MSO}_n(\mathrm{BDiff}^\delta(\Sigma_g)) \longrightarrow \mathcal{F}\Omega_{n+2,2}.$$

Using Theorem 0.4, we will determine the image of $e_2$ and $e_3$ up to torsions. More precisely,

**Theorem 0.18.** *For $g \geq 4$, the map*

$$e_3 : \mathrm{MSO}_3(\mathrm{BDiff}^\delta(\Sigma_g)) \longrightarrow \mathcal{F}\Omega_{5,2},$$

*is rationally surjective and for $g \geq 6$, it is rationally an isomorphism.*

*Remark* 0.19. To geometrically interpret the theorem, let $\mathcal{F}$ be any codimension 2 foliation on a manifold of dimension 5 and let $k\mathcal{F}$ denote a disjoint union of $k$ copies of $\mathcal{F}$. Then $k\mathcal{F}$, for some integer $k$, is foliated cobordant to a flat surface bundle of genus at most 4. And the injection for $g \geq 6$ means that if $\mathcal{F}$ is a flat surface bundle over a 3-manifold which bounds a codimension 2 foliation on a 6 manifold, then $k\mathcal{F}$ for some integer $k$ bounds a flat surface bundle over a 4-manifold where the genus of the fibers is at least 6.

*Remark* 0.20. This theorem might be compared to the result of Mizutani, Morita and Tsuboi in [MMT83] in which they proved that any codimension one foliation almost without holonomy is homologous to a disjoint union of flat circle bundles over tori.

**Theorem 0.21.** *Let $\phi$ be the map*

$$\phi \colon \mathcal{F}\Omega_{4,2} \otimes \mathbb{Q} \to \mathbb{Q},$$

*that sends a foliation $\mathcal{F}$ on a 4-manifold $M$ to the difference of the Pontryagin numbers $\int_M p_1(M) - p_1(\nu(\mathcal{F}))$. For $g \geq 3$, we have the following short exact sequence*

$$0 \to \mathrm{MSO}_2(\mathrm{BDiff}^\delta(\Sigma_g)) \otimes \mathbb{Q} \xrightarrow{e_2} \mathcal{F}\Omega_{4,2} \otimes \mathbb{Q} \xrightarrow{\phi} \mathbb{Q} \to 0.$$

*Remark* 0.22. Roughly speaking, up to torsion the only obstruction for a codimension 2 foliation $\mathcal{F}$ on a 4-manifold $M$ to be foliated cobordant to a flat surface bundle is $\int_M p_1(M) - p_1(\nu(\mathcal{F}))$.

We also prove that in low dimensions, we can change surface bundles up to cobordism to obtain a flat surface bundle; more precisely we prove

**Theorem 0.23.** *For $g > 5$, every surface bundle of genus $g$ over a 3-manifold is cobordant to a flat surface bundle.*

*Remark* 0.24. Using the perfectness of the identity component of $C^\infty$-diffeomorphisms, Morita and Kotschick in [KM05] proved that every surface bundle over a surface is foliated cobordant to a flat surface bundle.

0.5. **Outline.** This paper is organized as follows: In Section 1, we obtain a short proof of homological stability of discrete surface diffeomorphisms using a deep theorem of Mather and Thurston and a version of twisted stability of mapping class group due to Randal-Williams. In Section 2, we derive a Madsen-Weiss type theorem for discrete surface diffeomorphisms. In Section 3, we explore the consequences of having Madsen-Weiss type theorem for discrete surface diffeomorphisms in flat surface bundles and their characteristic classes.

**Acknowledgments.** It is my pleasure to thank my thesis advisor, Søren Galatius for his patience, encouragement and technical support during this project. Without his help and support, this paper would not have existed. I also would like to thank Oscar Randal-Williams for many sage remarks. I would like to thank Alexander Kupers for reading the first draft of this paper and Jonathan Bowden for pointing out to me an obstruction for a codimension 2 foliation on a 4-manifold to be foliated cobordant to a flat surface bundle. I am also very grateful to the referee whose many helpful comments improve the exposition of the paper. This work was partially supported by NSF grants DMS-1105058 and DMS-1405001.

## 1. HOMOLOGICAL STABILITY FROM FOLIATION THEORY

Our goal in this section is to show that the homological stability of discrete surface diffeomorphisms is implied by a "twisted" homological stability of mapping class groups developed in [RW16]. We use foliation theory to show that $\mathrm{BDiff}^\delta(\Sigma_{g,k}, \partial)$ is homology equivalent to a moduli space of a certain tangential structure in the sense of [RW16, Definition 1.1].

## 1.1. Stabilization maps.
In the introduction, we formulated the homological stability for discrete diffeomorphisms of $\Sigma_{g,1}$. Let us describe the stabilization maps for surfaces with any positive number of boundary components.

For a surface $\Sigma$ with boundary, we shall write $\mathrm{Diff}^\delta(\Sigma, \partial)$ to denote the discrete group of compactly supported orientation preserving diffeomorphisms of $\Sigma \backslash \partial \Sigma$ and if $\Sigma$ is a closed compact surface, $\mathrm{Diff}^\delta(\Sigma, \partial)$ will just mean all orientation preserving diffeomorphisms of $\Sigma$ equipped with discrete topology. Let $\Sigma \hookrightarrow \Sigma'$ be a subsurface in a collared surface $\Sigma'$ (i.e. with a choice of a collar neighborhood of the boundary) such that each of the boundary components of the subsurface either coincides with one of the boundary components of the bigger surface or entirely lies in its interior. If we extend diffeomorphisms of $\Sigma$ via the identity over the cobordism $K = \Sigma' \backslash \overline{\Sigma}$, we obtain a map

$$t: \mathrm{Diff}^\delta(\Sigma, \partial) \to \mathrm{Diff}^\delta(\Sigma', \partial).$$

Let $\Sigma_{g,k}$ be a fixed model for an orientable surface of genus $g$ and $k$ boundary components with a chosen collar neighborhood of the boundary. If we choose a diffeomorphism $f$ from $\Sigma$ to $\Sigma_{g,k}$ and a diffeomorphism $h$ from $\Sigma'$ to $\Sigma_{g',k'}$, we obtain a stabilization map

$$s_{f,h}(t): \mathrm{Diff}^\delta(\Sigma_{g,k}, \partial) \to \mathrm{Diff}^\delta(\Sigma_{g',k'}, \partial).$$

Different choices of $f$ and $h$ induce different stabilization maps, but one can show similar to [Nar14, Theorem 2.5] that if for some choices of $f$ and $h$, the map $s_{f,h}(t)$ induces a homology isomorphism in some homological degrees, then the stabilization map induces a homology isomorphism in the same homological degrees for all choices of $f$ and $h$. Therefore, we shall not write the dependence of the stabilization maps on choices of $f$ and $h$.

## 1.2. Homological stability of the moduli space of $\Gamma_2$-structures.
Because we are interested in homological stability of $\mathrm{Diff}^\delta(\Sigma_{g,k}, \partial)$, we may replace $\mathrm{BDiff}^\delta(\Sigma_{g,k}, \partial)$ by a homology equivalent space which is more convenient from the point of view of homotopy theory. To do so, we recall what we need from foliation theory.

### 1.2.1. Mather-Thurston's theory.
Recall from Definition 0.3 that $\mathrm{S\Gamma}_n$ is the groupoid of germs of orientation preserving diffeomorphisms of $\mathbb{R}^n$. Let $\mathrm{Mor}_{\mathrm{S\Gamma}_n}$ denote the space of morphisms in the topological groupoid $\mathrm{S\Gamma}_n$. To recall the topology on the space of morphisms, let $g \in \mathrm{Mor}_{\mathrm{S\Gamma}_n}$ be a germ sending $x$ to $y$. One can represent $g$ as a local diffeomorphism $\tilde{g}: U \to V$ where $U$ and $V$ are open sets containing $x$ and $y$ respectively. The set of germs of $\tilde{g}$ at all points in $U$ gives an open neighborhood of the germ $g$.

**Definition 1.1.** Let $X$ be a topological space. A 1-cocycle on $X$ with values in $\mathrm{S\Gamma}_n$ consists of an open cover $\{U_\alpha\}_I$ of $X$, and for any two indices $\alpha$ and $\beta$ in $I$, a continuous map $\gamma_{\alpha\beta}: U_\alpha \cap U_\beta \to \mathrm{Mor}_{\mathrm{S\Gamma}_n}$ satisfying cocycle condition for any $\alpha, \beta$ and $\delta$:

$$\gamma_{\alpha\beta}\gamma_{\beta\delta} = \gamma_{\alpha\delta} \text{ on } U_\alpha \cap U_\beta \cap U_\delta.$$

In particular, the left hand side is defined i.e. the source of the map $\gamma_{\alpha\beta}$ is the same as the target of $\gamma_{\beta\delta}$.

Two cocycles $c = \{U_\alpha, \gamma_{\alpha\beta}\}_I$ and $c' = \{U_{\alpha'}, \gamma_{\alpha'\beta'}\}_J$ are said to be equivalent if there exists a cocycle $c'' = \{U_{\alpha''}, \gamma_{\alpha''\beta''}\}_K$ such that $K = I \cup J$ and $c''$ restricts to $c$ on $\{U_\alpha\}_I$ and restricts to $c'$ on $\{U_{\alpha'}\}_J$.

**Definition 1.2.** An $S\Gamma_n$-structure on $X$ is an equivalence class of 1-cocycles with values in $\mathrm{Mor}_{S\Gamma_n}$ on $X$.

Note that a co-oriented foliation of codimension $n$ can be specified by a covering of $X$ by open sets $U_\alpha$, together with a submersion $f_\alpha$ from each open set $U_\alpha$ to $\mathbb{R}^n$, such that for each $\alpha$ and $\beta$ there is a map $g_{\alpha\beta}$ from $U_\alpha \cap U_\beta$ to local diffeomorphisms satisfying

$$f_\alpha(v) = g_{\alpha\beta}(u)\big(f_\beta(v)\big),$$

whenever $v$ is close enough to $u \in U_\alpha \cap U_\beta$. Then the covering $U_\alpha$ and the germs of $g_{\alpha\beta}$ defines a $S\Gamma_n$-structure on $X$. An advantage of Haefliger structures over foliations is that they are closed under pullbacks.

Two $S\Gamma_n$-structures on $X$, $c_0$ and $c_1$, are concordant if there exists a $S\Gamma_n$-structure on $X \times [0,1]$ so that the restriction of $c$ to $X \times \{i\}$ is $c_i$ for $i = 0, 1$. Homotopy classes of maps to the classifying space of the groupoid $S\Gamma_n$ classify $S\Gamma_n$-structures on $X$ up to concordance (for further details consult [Hae71]).

One can associate a foliated space to every $S\Gamma_n$-structure $c = \{U_\alpha, \gamma_{\alpha\beta}\}_I$ as follows. Note that for all $\alpha$, the cocycle condition implies that $\gamma_{\alpha\alpha}$ is the germ of the identity at some point in $\mathbb{R}^n$, hence $\gamma_{\alpha\alpha}$ induces a map from $U_\alpha$ to $\mathbb{R}^n$. Consider the space

$$\coprod_\alpha U_\alpha \times \mathbb{R}^n / \sim,$$

where the identification is given by $(x \in U_\alpha, y_\alpha) \sim (x \in U_\beta, y_\beta)$ if $y_\alpha = \gamma_{\alpha\beta} y_\beta$. We now consider $(x, \gamma_{\alpha\alpha}(x)) \in U_\alpha \times \mathbb{R}^n$, the graph of $\gamma_{\alpha\alpha}$. Let $E$ be the space obtained by the union of the neighborhoods of these graphs by the identification. The space $E$ is germinally well defined and the horizontal foliation on $U_\alpha \times \mathbb{R}^n$ induces a foliation on $E$. Hence, we obtain the data of a microbundle $X \xrightarrow{s} E \xrightarrow{p} X$ where $s$ is the section given by the graphs and $p$ is the projection to the first factor. The the foliation on $E$ is transverse to the fibers and its pullback to $X$ via the section $s$ is the Haefliger structure $c$. We call this microbundle, the foliated microbundle associated to $c$ (see [Tsu09, Section 4] for more details on foliated microbundles).

To state the Mather-Thurston theorem, we let $\overline{B S\Gamma}_n$ be the homotopy fiber of the natural map

$$\nu : B S\Gamma_n \to B\mathrm{GL}_n^+(\mathbb{R}),$$

which is induced by the map between groupoids that sends every germ to its derivative. By replacing spaces with homotopy equivalent spaces, we may assume that $\nu$ is a Serre fibration. For any orientable manifold $M$, we let $\overline{B\mathrm{Diff}(M, \partial)}$ be the homotopy fiber of the following map

$$B\mathrm{Diff}^\delta(M, \partial) \to B\mathrm{Diff}(M, \partial).$$

Let $\tau_M : M \to B\mathrm{GL}_n^+(\mathbb{R})$ be the map that classifies the tangent bundle of $M$. Thus, we obtain a Serre fibration $\tau_M^*(\nu)$ with the base $M$. The manifold structure on $M$ or the foliation by points on $M$ induces a homotopy class of maps $M \to B S\Gamma_n$. Let $s_0$ be one such map that is induced by the point foliation. Let also $\mathrm{Sect}_c(\tau_M^*(\nu))$ denote the compactly supported sections which differ only on a compact set from $s_0$.

**Theorem 1.3** (Mather-Thurston [Mat11])**.** *There exists a map*

$$f_M : \overline{B\mathrm{Diff}(M, \partial)} \to \mathrm{Sect}_c(\tau_M^*(\nu)),$$

*which induces an isomorphism in homology with integer coefficients.*

*Remark* 1.4. Roughly, the map $f_M$ is given by thinking of a diffeomorphism as a collection of germs at each point of $M$. Since the elements of $\overline{B\mathrm{Diff}(M, \partial)}$ can be thought as integrable sections of the fiber bundle $\tau_M^*(\nu)$, this theorem is very similar in spirit to Gromov h-principle type theorems.

*Remark* 1.5. Haefliger proved in [Hae71] that $\overline{\mathrm{B}\mathrm{S}\Gamma}_n$ is $n$-connected and he conjectured that it is $2n$-connected. Thurston could improve the connectivity of $\overline{\mathrm{B}\mathrm{S}\Gamma}_n$ by one. To explain his idea, note that Theorem 1.3 for an $n$-disk $D^n$ implies that the space $\overline{\mathrm{B}\mathrm{Diff}(D^n,\partial)}$ is homology isomorphic to the $n$-fold loop space $\Omega^n\overline{\mathrm{B}\mathrm{S}\Gamma}_n$. Given that Thurston ([Thu74]) also proved that the identity component of diffeomorphisms of manifolds is a perfect group (in fact he showed it is even simple), one can deduce that $H_1(\overline{\mathrm{B}\mathrm{Diff}(D^n,\partial)};\mathbb{Z}) = 0$. Therefore using Theorem 1.3 and the Hurewicz theorem, we obtain

$$H_1(\Omega^n\overline{\mathrm{B}\mathrm{S}\Gamma}_n;\mathbb{Z}) = H_{n+1}(\overline{\mathrm{B}\mathrm{S}\Gamma}_n;\mathbb{Z}) = \pi_{n+1}(\overline{\mathrm{B}\mathrm{S}\Gamma}_n) = 0.$$

As we shall see, the topological group $\mathrm{Diff}(M,\partial)$ acts on suitable models for $\overline{\mathrm{B}\mathrm{Diff}(M,\partial)}$ and $\mathrm{Sect}_c(\tau_M^*(\nu))$. Our goal is to show that the homotopy quotient of these actions are also homology equivalent. In order to achieve this goal, it is convenient to work with simplicial sets instead of topological spaces, and that we will explain how to define a map of simplicial sets modeling $f_M$, which is equivariant for an action of a simplicial group modeling $\mathrm{Diff}(M)$ (see [Nar14, Section 5.1] for a different model of the map $f_M$ which is equivariant). Henceforth, we substitute spaces with their singular simplicial complex.

1.2.2. *Construction of the map $f_M$.* Since $\overline{\mathrm{B}\mathrm{Diff}(M,\partial)}$ and $\mathrm{Sect}_c(\tau_M^*(\nu))$ classify certain geometric structures, it is more convenient to describe their singular simplicial complex geometrically. To do so, we need to recall few notions from [Ser79] and [Mat11].

**Definition 1.6.** We say a $\mathrm{S}\Gamma_n$-structure $c$ on the total space of the fiber bundle $E \to B$ is *transverse* to the fibers if its restriction to the fibers is a foliation. If the fiber bundle $E \to B$ is a smooth bundle, this is equivalent to the condition that $c$ is a smooth foliation and its leaves are transverse to the fibers.

**Definition 1.7.** Let $M$ be a smooth $n$-manifold, $X$ a topological space, and $c$ a $\mathrm{S}\Gamma_n$ on $X \times M$. We say $c$ is horizontal if $c$ is the inverse image of the differentiable structure of $M$ via the projection $X \times M \to M$. If $t \in X, x \in M$, we will say $c$ is locally horizontal at $(t,x)$ if there is an open neighborhood $U \times N$ of $(t,x)$ in $X \times M$ such that $c|_{U \times N}$ is horizontal. The support of $c$, denoted by $\mathrm{supp}(c)$, will mean the closure in $M$ of the set of $x \in M$ for which there is at least one $t \in X$ such that $c$ is not locally horizontal at $(t,x)$.

Since $\overline{\mathrm{B}\mathrm{Diff}(M,\partial)}$ is the homotopy fiber of the following map

$$\mathrm{B}\mathrm{Diff}^\delta(M,\partial) \to \mathrm{B}\mathrm{Diff}(M,\partial),$$

the $p$-simplices of the singular simiplicial complex, $S_\bullet(\overline{\mathrm{B}\mathrm{Diff}(M,\partial)})$, are uniquely given by $\mathrm{S}\Gamma_n$-structures on $\Delta^p \times M$ transverse to the fiber of the projection $\Delta^p \times M \to \Delta^p$ and have support in the interior of $M$.

The $p$-simplices of the simplicial group $S_\bullet(\mathrm{Diff}(M,\partial))$, namely the singular complex of $\mathrm{Diff}(M,\partial)$, can be described as the commutative diagrams

$$
\begin{array}{ccc}
\Delta^p \times M & \xrightarrow{\ \phi\ } & \Delta^p \times M \\
 & \searrow{\scriptstyle pr_1} \quad \swarrow{\scriptstyle pr_1} & \\
 & \Delta^p, &
\end{array}
$$

where $\phi$ is a diffeomorphism which is identity on $\Delta^p \times U$ where $U$ is a neighborhood of the boundary $\partial M$. We can pullback $S\Gamma_n$-structures on $\Delta^p \times M$ via $\phi$. Hence, we have an action of $S_\bullet(\text{Diff}(M, \partial))$ on $S_\bullet(\overline{\text{BDiff}(M, \partial)})$. Using the theorem of Milnor ([Mil57]), we know that $|S_\bullet(\overline{\text{BDiff}(M, \partial)})|$ is a model for $\overline{\text{BDiff}(M, \partial)}$, hence we obtain an action of the group $|S_\bullet(\text{Diff}(M, \partial))|$ which is weakly equivalent to $\text{Diff}(M, \partial)$ on $\overline{\text{BDiff}(M, \partial)}$. Therefore, the homotopy quotient [1]

$$(1.8) \qquad |S_\bullet(\overline{\text{BDiff}(M, \partial)})| /\!/ |S_\bullet(\text{Diff}(M, \partial))|$$

is weakly equivalent to $\text{BDiff}^\delta(M, \partial)$.

To describe the simplicial set $S_\bullet(\text{Sect}_c(\tau_M^*(\nu)))$ geometrically, we consider of the tangent bundle as the tangent microbundle (see [Mil64] for the definition of microbundles). Recall that the tangent microbundle of the manifold $M$ is the following data

$$M \xrightarrow{\Delta} M \times M \xrightarrow{pr_1} M,$$

which we denote by $\text{t}_M$. Milnor showed ([Mil64, Theorem 2.2]) that the underlying microbundle of the tangent bundle $TM$ is isomorphic to $\text{t}_M$. For an element $f \in \text{Diff}(M, \partial)$, we have an action of $f$ on $TM$ so that it acts by $f$ on $M$ and by the differential $df$ on the fiber of $TM$, the corresponding action on $\text{t}_M$ is acting by $f$ on the base $M$ and by $f \times f$ on $M \times M$.

Recall that every section in $\text{Sect}_c(\tau_M^*(\nu))$ is a lift of the tangent bundle in

$$
\begin{array}{ccc}
 & & \text{B}S\Gamma_n \\
 & \nearrow{\scriptstyle s} & \downarrow{\scriptstyle \nu} \\
M & \xrightarrow{\tau} & \text{BGL}_n^+(\mathbb{R})
\end{array}
$$

to $\text{B}S\Gamma_n$ which gives a map $s : M \to \text{B}S\Gamma_n$ and an isomorphism between $\text{t}_M$ and the underlying microbundle of $s^* \circ \nu^*(\gamma_n)$ where $\gamma_n$ is the tautological bundle on $\text{BGL}_n^+(\mathbb{R})$. This means that the graph of the $S\Gamma_n$-structure induced by $s$ is a foliated microbundle in the neighborhood of the diagonal $\Delta(M) \subset M \times M$ that is transverse to the fibers of $pr_1 : M \times M \to M$.

**Definition 1.9.** A germ of $S\Gamma_n$-structure $c$ on $\Delta^p \times M \times M$ at $\Delta^p \times \text{diag } M$ which is transverse to the fiber of the projection $id \times pr_1 : \Delta^p \times M \times M \to \Delta^p \times M$, is said to be horizontal at $x \in M$ if there exists a neighborhood $U$ around $x$ such that the restriction of the $S\Gamma_n$-structure to $\Delta^p \times U \times U$ is induced by the projection $\Delta^p \times U \times U \to U$ on the last factor. By support of $c$, we mean the set of $x \in M$, where $c$ is not horizontal. Note that $supp(c)$ is a closed subset.

Hence, $p$-simplices in the simplicial set $S_\bullet(\text{Sect}_c(\tau_M^*(\nu)))$ can be described as the germ of $S\Gamma_n$-structures on $\Delta^p \times M \times M$ at $\Delta^p \times \text{diag } M$ which are transverse to the fiber of the projection $id \times pr_1 : \Delta^p \times M \times M \to \Delta^p \times M$ and have compact support. This gives a model for the compactly supported sections $\text{Sect}_c(\tau_M^*(\nu))$ ([Mat11, Section 16]). Similar to the previous case, there is an obvious action of $S_\bullet(\text{Diff}(M, \partial))$ on $S_\bullet(\text{Sect}_c(\tau_M^*(\nu)))$.

**Construction 1.10.** Let $f_{M, \bullet} : S_\bullet(\overline{\text{BDiff}(M, \partial)}) \to S_\bullet(\text{Sect}_c(\tau_M^*(\nu)))$ be the simplicial map that sends a $p$-simplex $c$ in $S_\bullet(\overline{\text{BDiff}(M, \partial)})$, to the germ of the $S\Gamma_n$-structure induced by $(id \times pr_2)^*(c)$ at $\Delta^p \times \text{diag } M$, where $pr_2$ is the projection to the second factor. This map is obviously $S_\bullet(\text{Diff}(M, \partial))$-equivariant. Hence, using the Mather-Thurston theorem, the map $f_{M,\bullet}$ also induces homology isomorphism between homotopy quotients

$$|S_\bullet(\overline{\mathrm{BDiff}(M,\partial)})|/\!\!/|S_\bullet(\mathrm{Diff}(M,\partial))| \to |S_\bullet(\mathrm{Sect}_c(\tau_M^*(\nu)))|/\!\!/|S_\bullet(\mathrm{Diff}(M,\partial))|.$$

---

[1] For a topological group $G$ acting on a topological space $X$, the homotopy quotient is denoted by $X /\!/ G$ and is given by $X \times_G \text{E}G$ where $\text{E}G$ is a contractible space on which $G$ acts freely.

Therefore, Mather-Thurston's theorem imply that $\mathrm{BDiff}^\delta(M,\partial)$ is homology equivalent to $|S_\bullet(\mathrm{Sect}_c(\tau_M^*(\nu)))|/\!\!/|S_\bullet(\mathrm{Diff}(M,\partial))|$.

1.2.3. *Homological stability for tangential structures.* Recall from [RW16] that a *tangential structure* is a map $\theta : B \to \mathrm{BGL}_2(\mathbb{R})$ from a path-connected space $B$ to $\mathrm{BGL}_2(\mathbb{R})$. A $\theta$-structure on a surface $\Sigma$ is a bundle map $\mathrm{T}\Sigma \to \theta^*\gamma_2$ where $\gamma_2$ is the universal bundle over $\mathrm{BGL}_2(\mathbb{R})$. We denote the space of $\theta$-structures on a surface $\Sigma$ by $\mathrm{Bun}(\mathrm{T}\Sigma, \theta^*\gamma_2)$ and equip it with the compact-open topology. For a collared surface $\Sigma$ with a choice of collar $c : \partial\Sigma \times [0, 1) \to \Sigma$, we fix a boundary condition $\ell_{\partial\Sigma} : \epsilon^1 \oplus \mathrm{T}(\partial\Sigma) \to \theta^*\gamma_2$. And we define $\mathrm{Bun}_\partial(\mathrm{T}\Sigma, \theta^*\gamma_2; \ell_{\partial\Sigma})$ to be the space of bundle maps $\ell : \mathrm{T}\Sigma \to \theta^*\gamma_2$ such that $\ell_{\partial\Sigma} = Dc|_{\{0\}\times\partial\Sigma} \circ \ell|_{\partial\Sigma}$. Note that the group $\mathrm{Diff}(\Sigma, \partial)$ naturally acts on $\mathrm{Bun}_\partial(\mathrm{T}\Sigma, \theta^*\gamma_2; \ell_{\partial\Sigma})$. The moduli space of $\theta$-structures on surfaces of topological type $\Sigma$ with boundary condition $\ell_{\partial\Sigma}$ is the homotopy quotient of the action $\mathrm{Diff}(\Sigma, \partial)$ on $\mathrm{Bun}_\partial(\mathrm{T}\Sigma, \theta^*\gamma_2; \ell_{\partial\Sigma})$ and we denote it by

$$(1.11) \qquad \mathcal{M}^\theta(\Sigma; \ell_{\partial\Sigma}) := \mathrm{Bun}_\partial(\mathrm{T}\Sigma, \theta^*\gamma_2; \ell_{\partial\Sigma})/\!\!/\mathrm{Diff}(\Sigma, \partial).$$

If we do not mention the boundary condition $\ell_{\partial\Sigma}$, we mean the standard boundary condition on $\partial\Sigma$ in a sense of [RW16, Definition 4.1]. Henceforth, we consider $\nu$-structures where $\nu : \mathrm{B}S\Gamma_2 \to \mathrm{BGL}_2^+(\mathbb{R})$.

Recall that a foliated microbundle in the neighborhood of the diagonal of $\Sigma \times \Sigma$ which is transverse to the fibers of the projection $pr_1$, is a section in $\mathrm{Sect}_c(\tau_\Sigma^*(\nu))$. This then gives a bundle map from $T\Sigma$ to $\nu^*(\gamma_2)$. Therefore, there is a canonical map

$$\epsilon : |S_\bullet(\mathrm{Sect}_c(\tau_\Sigma^*(\nu)))| \xrightarrow{\simeq} \mathrm{Bun}_\partial(\mathrm{T}\Sigma, \nu^*\gamma_2; l_{\partial\Sigma}).$$

Using the fact that $\mathrm{Bun}_\partial(\mathrm{T}\Sigma, \gamma_2)$ is contractible ([GMTW09, Lemma 5.1]), one can show that there exists a homotopy inverse to the map $\epsilon$, hence it is a weak equivalence. There is an action of $|S_\bullet(\mathrm{Diff}(\Sigma, \partial))|$ on the left hand side and there is an action of $\mathrm{Diff}(\Sigma, \partial)$ on the right hand side and also there is a canonical map $\epsilon' : |S_\bullet(\mathrm{Diff}(\Sigma, \partial))| \to \mathrm{Diff}(\Sigma, \partial)$ which is weakly equivalent. The augmentation map $\epsilon$ is readily seen to be equivariant with respect to the map $\epsilon'$. Hence, using Serre spectral sequence we deduce that the induced map

$$(1.12) \qquad |S_\bullet(\mathrm{Sect}_c(\tau_\Sigma^*(\nu)))|/\!\!/|S_\bullet(\mathrm{Diff}(\Sigma, \partial))| \to \mathcal{M}^\nu(\Sigma),$$

is a homology isomorphism. As we saw in Construction 1.10, we have a map from $\mathrm{BDiff}^\delta(\Sigma, \partial)$ to $|S_\bullet(\mathrm{Sect}_c(\tau_\Sigma^*(\nu)))|/\!\!/|S_\bullet(\mathrm{Diff}(\Sigma, \partial))|$ which is homology equivalent. For the future reference, we record this fact as a lemma.

**Lemma 1.13.** *There is a map from $\mathrm{BDiff}^\delta(\Sigma, \partial)$ to $\mathcal{M}^\nu(\Sigma)$ that induces an isomorphism on homology.*

Hence Lemma 1.13 reduces the proof of Theorem 0.1 to the homological stability of the moduli space of $\nu$-structures, which then follows from a general theorem due to Randal-Williams [RW16, Theorem 7.1] about homological stability of moduli spaces of $\theta$-structures satisfying certain properties. Qualitatively, he proved that if the connected components of the $\mathcal{M}^\nu(\Sigma)$ stabilizes with respect to the genus of the surface $\Sigma$, the moduli space $\mathcal{M}^\nu(\Sigma)$ exhibits homological stability in a certain range.

To show the stability of connected components, consider the following exact sequence of homotopy groups

$$\pi_1(\mathrm{BDiff}(\Sigma,\partial)) \to \pi_0(\mathrm{Sect}_c(\tau_\Sigma^*(\nu))) \to \pi_0(\mathcal{M}^\nu(\Sigma)) \to \pi_0(\mathrm{BDiff}(\Sigma,\partial)).$$

The classifying space $\mathrm{BDiff}(\Sigma,\partial)$ is path-connected and since by Remark 1.5 the space $\overline{\mathrm{B}\mathrm{S}\Gamma_2}$ is at least 3-connected, the section space $\mathrm{Sect}_c(\tau_\Sigma^*(\nu))$ is also path-connected. Hence, $\pi_0(\mathcal{M}^\nu(\Sigma))$ is trivial.

To find a stability range, Randal-Williams defined a notion of $k$-triviality [RW16, Definition 6.2] and proved that if a $\theta$-structure stabilizes at genus $h$, then it would be $(2h+1)$-trivial. Since $\nu$-structure stabilizes at genus 0, by [RW16, Theorem 7.1] the stability range for stabilization maps is the same as the stability range for the orientation structure $\mathrm{BSO}(2) \to \mathrm{BO}(2)$. Thus, we have

**Theorem 1.14.** *Let $\Sigma$ be a collared surface and let $\iota : \Sigma \hookrightarrow \Sigma'$ be an embedding of $\Sigma$ into a surface $\Sigma'$ which may not have boundary. As we explained in Section 1.1, this embedding induces a map*

$$H_*(\mathrm{BDiff}^\delta(\Sigma,\partial);\mathbb{Z}) \to H_*(\mathrm{BDiff}^\delta(\Sigma',\partial);\mathbb{Z})$$

*which is an isomorphism as long as $* \le (2g(\Sigma) - 2)/3$ and epimorphism provided that $* \le 2g(\Sigma)/3$.*

*Remark* 1.15. Using the same idea and the theorem of McDuff [McD83] about volume preserving diffeomorphisms, we could show that the discrete group of symplectomorphisms $\mathrm{Symp}^\delta(\Sigma,\partial)$ that are supported away from the boundary also exhibit homological stability. We pursue the study the group of symplectomorphisms in a different paper ([Nar]).

## 2. Stable homology of surface diffeomorphisms made discrete

Given that we established the relation between $\mathrm{BDiff}^\delta(\Sigma,\partial)$ for a collared surface $\Sigma$ and the moduli space of $\mathrm{S}\Gamma_2$-structures on surfaces of the topological type $\Sigma$ in Lemma 1.13, we can use the machinery developed in [GMTW09], [GRW10] to study the stable homology of the moduli space of tangential structures.

### 2.1. Cobordism category with $\mathrm{S}\Gamma_2$-structure.
Recall the definition of the cobordism category equipped with $\theta$-structures from [GMTW09, Definition 5.2];

**Definition 2.1.** Let $\mathcal{C}_\theta$ be the topological category whose space of objects is given by the pairs of a real number $a$ and 1-dimensional closed submanifolds $M$ such that $(a, M) \subset \mathbb{R} \times \mathbb{R}^\infty$ and whose space of morphisms from $(a, M)$ to $(a', M')$ for $a < a'$ is given by a cobordism $\Sigma$ so that $\Sigma \subset [a, a'] \times \mathbb{R}^\infty$ is a surface equipped with a $\theta$-structure and is *collared* near the boundary which means that it coincides with $[a, a'] \times M$ near $\{a\} \times \mathbb{R}^\infty$ and with $[a, a'] \times M'$ near $\{a'\} \times \mathbb{R}^\infty$. For the careful treatment of how this category is enriched over topological spaces consult [GRW10, Section 2]. We shall write $\mathcal{C}_+$ for the cobordism category with the orientation structure.

**Definition 2.2.** For the map $\nu : \mathrm{BS}\Gamma_2 \to \mathrm{BGL}_2^+(\mathbb{R}) = \widetilde{\mathrm{Gr}}_2(\mathbb{R}^\infty)$ where $\widetilde{\mathrm{Gr}}_2(\mathbb{R}^\infty)$ is the oriented Grassmannian of two planes in $\mathbb{R}^\infty$, we can associate a Thom spectrum $\mathrm{MT}\nu$ as follows: First let $\mathrm{BS}\Gamma_2(\mathbb{R}^n) := \nu^{-1}(\widetilde{\mathrm{Gr}}_2(\mathbb{R}^n))$ where it sits in the pullback diagram

$$\begin{array}{ccc} \mathrm{BS}\Gamma_2(\mathbb{R}^n) & \longrightarrow & \mathrm{BS}\Gamma_2 \\ {\scriptstyle \nu_n} \downarrow & & \downarrow {\scriptstyle \nu} \\ \widetilde{\mathrm{Gr}}_2(\mathbb{R}^n) & \longrightarrow & \mathrm{BGL}_2^+(\mathbb{R}). \end{array}$$

Let $U_n$ be the orthogonal complement of the tautological 2-plane bundle over $\widetilde{\mathrm{Gr}}_2(\mathbb{R}^n)$. Then the $n$-th space of the spectrum $\mathrm{MT}\nu$ is the Thom space of the pullback bundle $\nu_n^*(U_n)$.

The main theorem [GMTW09] implies that there exists a weak equivalence

$$\mathrm{B}\mathcal{C}_\nu \xrightarrow{\simeq} \Omega^{\infty-1}\mathrm{MT}\nu$$

which is induced by a functor from the category $\mathcal{C}_\nu$ to the category $\mathrm{Path}(\Omega^{\infty-1}\mathrm{MT}\nu)$, whose objects are points in $\Omega^{\infty-1}\mathrm{MT}\nu$ and whose morphisms are continuous paths. We shall briefly recall below how this functor is constructed and refer the reader to [MT01, Section 2] for further details.

## 2.2. The map $\alpha$ in Theorem 0.4.

A morphism $\Sigma$ in the category $\mathcal{C}_\nu$ is a surface with a collared boundary that is embdded in $[a, a'] \times \mathbb{R}^{n-1}$ for some $n$. We say that this morphism is *fatly embedded* if the canonical map from the normal bundle $N\Sigma$ to $\mathbb{R}^n$ restricts to an embedding of the unit disk bundle into $[a, a'] \times \mathbb{R}^{n-1}$. In the definition of the cobordism category, one can consider only fatly embedded morphisms without changing the homotopy type of the realization of the category. Thus the Pontryagin-Thom construction give a map from $[a, a']_+ \wedge S^{n-1}$ to the Thom space of $N\Sigma$. Since the surface $\Sigma$ is equipped with a $\nu$-structure, the Gauss map $\Sigma \to \widetilde{\mathrm{Gr}}_2(\mathbb{R}^n)$ that classifies the tangent bundle can be lifted to $\mathrm{B}\mathrm{S}\Gamma_2(\mathbb{R}^n)$. Therefore, we have the pullback diagram

$$
\begin{array}{ccc}
N\Sigma & \longrightarrow & \nu_n^*(U_n) \\
\downarrow & & \downarrow \\
\Sigma & \longrightarrow & \mathrm{B}\mathrm{S}\Gamma_2(\mathbb{R}^n).
\end{array}
$$

Hence, one obtains the map

$$[a, a']_+ \wedge S^{n-1} \to \mathrm{Th}(N\Sigma) \to \mathrm{Th}(\nu_n^*(U_n)).$$

By the adjointness, we obtain a path

$$[a, a'] \to \Omega^{n-1}\mathrm{Th}(\nu_n^*(U_n)) \subset \Omega^{\infty-1}\mathrm{MT}\nu.$$

This construction gives rise to a functor from the modification of $\mathcal{C}_\nu$ to the path category $\mathrm{Path}(\Omega^{\infty-1}\mathrm{MT}\nu)$. Since the modification of $\mathcal{C}_\nu$ does not change its homotopy type and the geometric realization of $\mathrm{Path}(\Omega^{\infty-1}\mathrm{MT}\nu)$ has the same homotopy type as $\Omega^{\infty-1}\mathrm{MT}\nu$, the functor induces a well-defined map up to homotopy between geometric realizations

$$\mathrm{B}\mathcal{C}_\nu \to \Omega^{\infty-1}\mathrm{MT}\nu.$$

One can choose a certain model for the homotopy quotient in 1.11 (see [GMTW09, Section 5]) so that the space $\mathcal{M}^\nu(\Sigma)$ becomes a subspace of the morphism space in $\mathcal{C}_\nu$. Therefore, we obtain a natural map

$$(2.3) \qquad \mathcal{M}^\nu(\Sigma) \to \Omega\mathrm{B}\mathcal{C}_\nu \to \Omega^\infty\mathrm{MT}\nu.$$

Note that the map $\nu : \mathrm{B}\mathrm{S}\Gamma_2 \to \mathrm{B}\mathrm{GL}_2^+(\mathbb{R})$ induces a functor $\mathcal{C}_\nu \to \mathcal{C}_+$, hence by the naturality of the above constructions, we have the homotopy commutative diagram

$$
\begin{array}{ccccc}
\mathcal{M}^\nu(\Sigma) & \longrightarrow & \Omega\mathrm{B}\mathcal{C}_\nu & \longrightarrow & \Omega^\infty\mathrm{MT}\nu \\
\downarrow & & \downarrow & & \downarrow \\
\mathcal{M}^+(\Sigma) & \longrightarrow & \Omega\mathrm{B}\mathcal{C}_+ & \longrightarrow & \Omega^\infty\mathrm{MTSO}(2).
\end{array}
$$

Recall that the space $\mathcal{M}^+(\Sigma)$ is a model for $\mathrm{BDiff}(\Sigma, \partial)$ and the space

$$|S_\bullet(\overline{\mathrm{BDiff}(M, \partial)})|/\!\!/|S_\bullet(\mathrm{Diff}(M, \partial))|$$

is a model for $\mathrm{BDiff}^\delta(\Sigma, \partial)$. Hence, we have the following homotopy commutative diagram

$$(2.4) \qquad \begin{array}{ccc}
|S_\bullet(\overline{\mathrm{BDiff}(M, \partial)})|/\!\!/|S_\bullet(\mathrm{Diff}(M, \partial))| & \xrightarrow{\ g\ } & \mathcal{M}^\nu(\Sigma) \\
\downarrow & & \downarrow \\
\mathrm{B}|S_\bullet(\mathrm{Diff}(M, \partial))| & \xrightarrow{\ \simeq\ } & \mathcal{M}^+(\Sigma),
\end{array}$$

where the bottom horizontal map is a weak equivalence by Milnor's theorem [Mil57] and the top horizontal map $g$ is given by the composition of the map in Construction 1.10 and 1.12. The map $\alpha$ now is given by the composition of the map $g: \mathrm{BDiff}^\delta(\Sigma) \to \mathcal{M}^\nu(\Sigma)$ in 2.4 and the maps in 2.3. Hence, we obtain a homotopy commutative diagram

$$(2.5) \qquad \begin{array}{ccc}
\mathrm{BDiff}^\delta(\Sigma, \partial) & \xrightarrow{\ \alpha\ } & \Omega^\infty \mathrm{MT}\nu \\
\downarrow & & \downarrow \\
\mathrm{BDiff}(\Sigma, \partial) & \longrightarrow & \Omega^\infty \mathrm{MTSO}(2).
\end{array}$$

**Theorem 2.6.** *In the diagram 2.5 the horizontal maps, in the stable range of Theorem 1.14, induce homology isomorphisms onto the connected components that they hit.*

*Remark* 2.7. The volume preserving case reproduces [KM07, Theorem 4] and more which we will pursue elsewhere ([Nar]).

*Sketch of the proof of theorem 2.6.* The fact that the bottom horizontal map in the stable range induces an isomorphisms on homology is the celebrated Madsen-Weiss theorem ([MW07], [GMTW09, Theorem 7.2]). Hence, we only sketch the proof for the similar statement for the map $\alpha$. We replace $\mathrm{BDiff}^\delta(\Sigma_{g,k}, \partial)$ by the homology equivalent space $\mathcal{M}^\nu(\Sigma_{g,k})$. Recall that the main theorem of [GMTW09] implies that the geometric realization of $\mathcal{C}_\nu$ is weakly homotopy equivalent to $\Omega^{\infty-1}\mathrm{MT}\nu$. Therefore from the above discussion, we only need to prove that the map

$$(2.8) \qquad\qquad \mathcal{M}^\nu(\Sigma_{g,k}) \to \Omega\mathrm{B}\mathcal{C}_\nu,$$

in the stable range induces an isomorphism on homology. As we shall briefly explain, this follows from applying the argument in [GMTW09, Section 7] to the category $\mathcal{C}_\nu$. Following Tillmann [Til97], we need to consider a smaller category which is called *the positive boundary* subcategory $\mathcal{C}_{\nu,\partial} \subset \mathcal{C}_\nu$ whose space of objects is the same as $\mathcal{C}_\nu$ and whose space of morphisms from $M_0$ to $M_1$ consists of those pairs $(\Sigma, t) \in \mathcal{C}_\nu$ where $\pi_0(M_1) \to \pi_0(\Sigma)$ is surjective. By [GMTW09, Theorem 6.1] the inclusion of $\mathcal{C}_{\nu,\partial}$ in $\mathcal{C}_\nu$ induces a map between geometric realizations

$$(2.9) \qquad\qquad \mathrm{B}\mathcal{C}_{\nu,\partial} \xrightarrow{\ \simeq\ } \mathrm{B}\mathcal{C}_\nu,$$

which is a weak equivalence. There is nothing special about the tangential structure $\nu$ in 2.9. Now given the homological stability in Theorem 1.14, the standard group completion argument (see [GMTW09, Proposition 7.1] ) implies that the first map in

$$\mathbb{Z} \times \mathcal{M}^\nu(\Sigma_{\infty,k}) \xrightarrow{H_*-\mathrm{iso}} \Omega\mathrm{B}\mathcal{C}_{\nu,\partial} \xrightarrow{\ \simeq\ } \Omega^\infty\mathrm{MT}\nu,$$

induces an isomorphism on homology, and hence the map in 2.8 induces a homology isomorphism in the stable range.                                                      □

*Remark* 2.10. There is a more direct description of the diagram 2.5 without invoking the cobordism category. To briefly explain this alternative description, let $\Sigma \to E \xrightarrow{\pi} M$ be a surface bundle and let $T\pi$ denote the vertical tangent bundle. For paracompact base $M$, one can find a fiberwise embedding $f : E \hookrightarrow M \times \mathbb{R}^N$ with a tubular neighborhood. Collapsing the complement of the tubular neighborhood to the base point and identifying the tubular neighborhood with the open disk bundle of the fiberwise normal bundle $Nf$, gives the map $\Sigma^N(M_+) \to \mathrm{Th}(Nf)$. Stably it gives *a pretransfer map* well-defined up to homotopy

$$\mathrm{pretrf}_\pi : \Sigma^\infty(M_+) \to \mathbf{Th}(-T\pi),$$

where $\mathbf{Th}(-T\pi)$ is the Thom spectrum of the virtual bundle $-T\pi$. Let $\theta : B \to \mathrm{BGL}_2^+(\mathbb{R})$ be a tangential structure. Recall that $\gamma_2$ is the tautological bundle over $\mathrm{BGL}_2^+(\mathbb{R})$. If the vertical tangent bundle is equipped with a $\theta$-structure i.e.if the map $E \to \mathrm{BGL}_2^+(\mathbb{R})$ classifying $T\pi$, has a choice of lift to $B$, then we obtain a well-defined map up to homotopy

$$(2.11) \qquad \Sigma^\infty(M_+) \xrightarrow{\mathrm{pretrf}_\pi} \mathbf{Th}(-T\pi) \to \mathbf{Th}(-\theta^*(\gamma_2)).$$

Now consider the following pullback diagram

$$\begin{array}{ccc}
\Sigma /\!\!/ \mathrm{Diff}^\delta(\Sigma, \partial) & \longrightarrow & \Sigma /\!\!/ \mathrm{Diff}(\Sigma, \partial) \\
{\scriptstyle \pi'} \downarrow & & \downarrow {\scriptstyle \pi} \\
\mathrm{BDiff}^\delta(\Sigma, \partial) & \longrightarrow & \mathrm{BDiff}(\Sigma, \partial).
\end{array}$$

Note that $T\pi'$ is the pullback of $T\pi$ and it has a $\nu$-structure, therefore by the naturality of the pretransfer map, we have a commutative diagram of spectra

$$\begin{array}{ccc}
\Sigma^\infty(\mathrm{BDiff}^\delta(\Sigma, \partial)_+) & \to & \mathbf{Th}(-\nu^*(\gamma_2)) \\
\downarrow & & \downarrow \\
\Sigma^\infty(\mathrm{BDiff}(\Sigma, \partial)_+) & \longrightarrow & \mathbf{Th}(-\gamma_2),
\end{array}$$

which gives a homotopy commutative diagram of spaces

$$\begin{array}{ccc}
\mathrm{BDiff}^\delta(\Sigma, \partial) & \longrightarrow & \Omega^\infty \mathrm{MT}\nu \\
\downarrow & & \downarrow \\
\mathrm{BDiff}(\Sigma, \partial) & \longrightarrow & \Omega^\infty \mathrm{MTSO}(2).
\end{array}$$

The fact this diagram is the same as the diagram 2.5 follows from the standard fact that the construction using pretransfer and the construction via the cobordism category are homotopic (see [MT01, Section 2] for more details).

## 2.3. Comparison of $\mathrm{BDiff}^\delta(\Sigma_{g,k}, \partial)$ and $\mathrm{BDiff}(\Sigma_{g,k}, \partial)$ in the stable range.
Let $\mathrm{Diff}^\delta(\Sigma_{\infty,k}, \partial)$ denote the colimit of the groups $\mathrm{Diff}^\delta(\Sigma_{g,k}, \partial)$ as $g$ varies using the stabilization map between them. Thus by taking the colimit of the diagram

2.5, we have a homotopy commutative diagram

$$
\begin{array}{ccc}
\mathrm{BDiff}^{\delta}(\Sigma_{\infty,k}, \partial) & \xrightarrow{\iota} & \mathrm{BDiff}(\Sigma_{\infty,k}, \partial) \\
\downarrow & & \downarrow \\
\Omega_0^{\infty}\mathrm{MT}\nu & \xrightarrow{\Omega^{\infty}\boldsymbol{\nu}} & \Omega_0^{\infty}\mathrm{MTSO}(2),
\end{array}
$$

where the right vertical map is a homology isomorphism by the Madsen-Weiss theorem and the left vertical map is also a homology isomorphism by theorem 2.6. For a prime $p$, let

$$
\iota_*^p : H_*(\mathrm{BDiff}^{\delta}(\Sigma_{\infty,k}, \partial); \mathbb{F}_p) \to H_*(\mathrm{BDiff}(\Sigma_{\infty,k}, \partial); \mathbb{F}_p)
$$

be the map induced by $\iota$. To study the map $\iota_*^p$, instead we study the map

$$
\Omega^{\infty}\boldsymbol{\nu} : \Omega_0^{\infty}\mathrm{MT}\nu \to \Omega_0^{\infty}\mathrm{MTSO}(2),
$$

between infinite loop spaces after $p$-completion (see [MP11] for a definition of $p$-completion).

Our third main theorem is the following splitting theorem after $p$-completion.

**Theorem 2.12.** *For all prime $p$, after $p$-completion, the map $\Omega^{\infty}\boldsymbol{\nu}$ admits a section i.e. it is split surjection after $p$-adic completion.*

**Corollary 2.13.** *The map $\iota_*^p$ is a split surjection for all prime $p$.*

*Proof of Theorem 2.12.* To show that $\Omega^{\infty}\boldsymbol{\nu}$ has a section after $p$-completion, it is sufficient to prove that on the spectrum level, i.e. as we shall see, it is enough to prove that the map

$$
(2.14) \qquad\qquad \boldsymbol{\nu} : \mathrm{MT}\nu \to \mathrm{MTSO}(2)
$$

is a split surjection after $p$-completion of spectra. The reason is, in general, if $A$ is a spectrum and $A_p^{\wedge}$ is its $p$-completion, then $\Omega_0^{\infty}(A_p^{\wedge})$ is a $p$-completed space and it receives a map from $(\Omega_0^{\infty}A)_p^{\wedge}$ which is weakly equivalent. The map is induced by

$$
\Omega_0^{\infty}A \to \Omega_0^{\infty}(A_p^{\wedge})
$$

which factors through $(\Omega_0^{\infty}A)_p^{\wedge}$ by universal property of $p$-completion. Suppose the map $\boldsymbol{\nu}_p^{\wedge}$ which is induced by $\boldsymbol{\nu}$ between $p$-completions of the Thom spectra, has a section denoted by $\boldsymbol{s}$. Consider the diagram

$$
\begin{array}{ccc}
\Omega_0^{\infty}(\mathrm{MTSO}(2))_p^{\wedge} & \xrightarrow{\Omega^{\infty}\boldsymbol{s}} & \Omega_0^{\infty}(\mathrm{MT}\nu)_p^{\wedge} \\
\simeq \uparrow & & \uparrow \simeq \\
(\Omega_0^{\infty}\mathrm{MTSO}(2))_p^{\wedge} & \dashrightarrow & (\Omega_0^{\infty}\mathrm{MT}\nu)_p^{\wedge}
\end{array}
$$

where the vertical maps are induced by the universal property of the $p$-completion. Hence, $\Omega^{\infty}\boldsymbol{s}$ induces a section for $\Omega^{\infty}\boldsymbol{\nu}$ after $p$-completion. We are left to prove the map of spectra 2.14 is a split surjection after $p$-completion.

Recall that the map

$$
\nu : B\mathrm{S}\Gamma_2 \to \mathrm{BGL}_2^+(\mathbb{R}) \simeq \mathrm{BSO}(2)
$$

that classifies the normal bundle of the codimension 2 Haefliger structures, is induced by a map between groupoids $\bar{\nu} : \mathrm{S}\Gamma_2 \to \mathrm{GL}_2^+(\mathbb{R})$ where $\bar{\nu}$ sends a germ to its derivative.

The key point is that there is an obvious map between groupoids

$$\bar{\alpha} : \mathrm{SO}(2)^{\delta} \to \mathrm{S}\Gamma_2$$

that sends every rotation to its germ as a diffeomorphism at the origin and make the diagram

(2.15)

$$\mathrm{SO}(2)^{\delta} \xrightarrow{\bar{\alpha}} \mathrm{S}\Gamma_2$$
$$\bar{\beta} \searrow \quad \downarrow \bar{\nu}$$
$$\mathrm{GL}_2(\mathbb{R})$$

commute where $\bar{\beta}$ is the composition $\mathrm{SO}(2)^{\delta} \to \mathrm{SO}(2) \hookrightarrow \mathrm{GL}_2(\mathbb{R})$. The map $\bar{\alpha}$ induces a map between classifying spaces of the groupoids which we denote by $\alpha$. The commutativity of 2.15 implies that the composition $\nu \circ \alpha$ is homotopic to the identity.

The map $\nu \circ \alpha$ gives a tangential structure and let $\mathrm{MT}(\nu \circ \alpha)$ denote the Thom spectrum of the virtual bundle $(\nu \circ \alpha)^*(-\gamma_2)$ over $\mathrm{BSO}(2)^{\delta}$. Consider the following maps between Thom spectra

$$\mathrm{MT}(\nu \circ \alpha) \xrightarrow{\boldsymbol{\alpha}} \mathrm{MT}\nu \xrightarrow{\boldsymbol{\nu}} \mathrm{MTSO}(2).$$

If we show that the map of spectra $\boldsymbol{\nu} \circ \boldsymbol{\alpha}$ is a split surjection after $p$-completion, we are done. We prove that $(\boldsymbol{\nu} \circ \boldsymbol{\alpha})^*$ is an isomorphism on mod $p$ cohomology, hence $\boldsymbol{\nu} \circ \boldsymbol{\alpha}$ actually induces a weak equivalence after $p$-completion. Let $\overline{\mathrm{BS}^1}$ denote the homotopy fiber of the map

$$\mathrm{B}(\mathrm{S}^1)^{\delta} \to \mathrm{BS}^1.$$

It is a special case of [Mil83, Lemma 3] that the space $\overline{\mathrm{BS}^1}$ has a mod $p$ homology of a point. Thus, we deduce that $H^*(\mathrm{B}(\mathrm{S}^1)^{\delta}; \mathbb{F}_p) = H^*(\mathrm{BS}^1; \mathbb{F}_p)$. Hence using Thom isomorphism, it follows that

$$(\boldsymbol{\nu} \circ \boldsymbol{\alpha})^* : H^*(\mathrm{MTSO}(2); \mathbb{F}_p) \to H^*(\mathrm{MT}(\nu \circ \alpha); \mathbb{F}_p)$$

is an isomorphism. □

## 3. Applications to flat surface bundles

In this section, we explore the consequences of Theorem 2.6 and Theorem 2.12 for flat surface bundles.

### 3.1. The spectrum $\mathrm{MT}\nu$ and a fiber sequence of infinite loop spaces.
By studing the homotopy groups of $\Omega_0^{\infty}\mathrm{MT}\nu$, we prove Theorem 0.18 and also we find a new description for $H_2(\mathrm{BDiff}^{\delta}(\Sigma_g); \mathbb{Z})$.

Let $\mathrm{MTSO}(n)$ denote the Madsen-Tillmann spectrum for the orientation structure $\mathrm{SO}(n) \to \mathrm{O}(n)$. There exists a cofiber sequence of spectra and a fiber sequence of infinite loop space (see [GMTW09, Proposition 3.1]) as follows

(3.1) $$\mathrm{MTSO}(n) \to \Sigma^{\infty}(\mathrm{BSO}(n)_+) \to \mathrm{MTSO}(n-1)$$

(3.2) $$\Omega^{\infty}\mathrm{MTSO}(n) \to \Omega^{\infty}\Sigma^{\infty}(\mathrm{BSO}(n)_+) \to \Omega^{\infty}\mathrm{MTSO}(n-1)$$

where here + means a disjoint base point. For $n = 2$, the fiber sequence of the infinite loop space plays an important role in computing mod $p$ homology of the mapping class group [Gal04]. In dimension 2, we prove that there exists a similar fiber sequence for $\Omega^{\infty}\mathrm{MT}\nu$.

**Theorem 3.3.** *There is a homotopy fibration sequence*

(3.4) $$\Omega^{\infty}\mathrm{MT}\nu \to \Omega^{\infty}\Sigma^{\infty}((\mathrm{B}\mathrm{S}\Gamma_2)_+) \to \Omega^{\infty}\Sigma^{\infty-1}(\overline{\mathrm{B}\mathrm{S}\Gamma_2}_+)$$

*Proof.* Let $\eta$ and $\xi$ be two vector bundles over a topological space $X$, we have the following general cofiber sequence of Thom spaces (see [Tei92, Lemma 4.3.1])

$$(3.5) \qquad \mathrm{Th}(p^*(\xi)) \to \mathrm{Th}(\xi) \to \mathrm{Th}(\xi \oplus \eta)$$

where $p: S(\eta) \to X$ is the sphere bundle of $\eta$.

Let $X$ be $\mathrm{B}\Gamma_2$ and $\xi$ and $\eta$ be the virtual bundles $\nu^*(-\gamma_2)$ and $\nu^*(\gamma_2)$ respectively. Note that

$$S^1 \to \overline{\mathrm{B}\Gamma_2} \xrightarrow{p} \mathrm{B}\Gamma_2$$

is the sphere bundle of $\nu^*(\gamma_2)$, hence $p^*(\nu^*(-\gamma_2))$ is a trivial bundle over $\overline{\mathrm{B}\Gamma_2}$. Using 3.5, we obtain the following cofiber sequence of the spectra

$$\Sigma^{\infty-2}(\overline{\mathrm{B}\Gamma_{2+}}) \to \mathrm{MT}\nu \to \Sigma^\infty((\mathrm{B}\Gamma_2)_+) \to \Sigma^{\infty-1}(\overline{\mathrm{B}\Gamma_{2+}})$$

Applying $\Omega^\infty$, we obtain the associated homotopy fibration sequence of infinite loop spaces.  $\square$

**Lemma 3.6.** *We have*

$$\pi_1(\Omega_0^\infty \mathrm{MT}\nu) = 0$$
$$0 \to \pi_4(\overline{\mathrm{B}\Gamma_2}) \to \pi_2(\Omega_0^\infty \mathrm{MT}\nu) \to \pi_2(\Omega_0^\infty \mathrm{MTSO}(2)) \to 0$$
$$\pi_3(\Omega_0^\infty \mathrm{MT}\nu) \twoheadrightarrow \pi_3(\Omega_0^\infty \mathrm{MTSO}(2)).$$

*Proof.* Recall that for $g \geq 3$, the group $\mathrm{Diff}^\delta(\Sigma_g)$ is a perfect group because the identity component $\mathrm{Diff}_0^\delta(\Sigma_g)$ is a simple group ([Thu74]) and the mapping class group $\pi_0(\mathrm{Diff}(\Sigma_g))$ is perfect for $g \geq 3$ ([Pow78]). Hence $H_1(\mathrm{BDiff}^\delta(\Sigma_g); \mathbb{Z}) = 0$ for $g \geq 3$. On the other hand, Theorem 2.6 implies that $H_1(\mathrm{BDiff}^\delta(\Sigma_g); \mathbb{Z}) = H_1(\Omega_0^\infty \mathrm{MT}\nu; \mathbb{Z}) = \pi_1(\Omega_0^\infty \mathrm{MT}\nu)$, where the second equality holds because $\Omega_0^\infty \mathrm{MT}\nu$ is an H-space. Therefore we have

$$\pi_1(\Omega_0^\infty \mathrm{MT}\nu) = 0.$$

Also recall that for stable homotopy groups, we have $\pi_i^s(X_+) = \pi_i^s(X) \oplus \pi_i^s$, $\pi_4^s = \pi_5^s = 0$, $\pi_3^s = \mathbb{Z}/24$ and $\pi_2^s = \mathbb{Z}/2$. Given these facts, the long exact sequence of the homotopy groups of the fibration 3.4 is as follows

$$\cdots \longrightarrow \pi_5^s(\overline{\mathrm{B}\Gamma_2}) \longrightarrow \pi_3(\Omega_0^\infty \mathrm{MT}\nu) \longrightarrow \pi_3^s(\mathrm{B}\Gamma_2) \oplus \mathbb{Z}/24$$
$$\hookrightarrow \pi_4^s(\overline{\mathrm{B}\Gamma_2}) \longrightarrow \pi_2(\Omega_0^\infty \mathrm{MT}\nu) \longrightarrow \pi_2^s(\mathrm{B}\Gamma_2) \oplus \mathbb{Z}/2 \longrightarrow \pi_3^s(\overline{\mathrm{B}\Gamma_{2+}}).$$

By Remark 1.5, we know that the map $\nu: \mathrm{B}\Gamma_2 \to \mathrm{BSO}(2)$ is at least 4-connected. Hence, the long exact sequence of the homotopy groups of the fibrations 3.4 becomes

$$(3.7) \qquad \cdots \longrightarrow \pi_5^s(\overline{\mathrm{B}\Gamma_2}) \longrightarrow \pi_3(\Omega_0^\infty \mathrm{MT}\nu) \longrightarrow \mathbb{Z}/24$$
$$\hookrightarrow \pi_4^s(\overline{\mathrm{B}\Gamma_2}) \longrightarrow \pi_2(\Omega_0^\infty \mathrm{MT}\nu) \longrightarrow \pi_2^s(\mathbb{CP}^\infty) \oplus \mathbb{Z}/2 \xrightarrow[e]{d} \mathbb{Z}/24,$$

and the long exact sequence of homotopy groups of the fibration in 3.2 is as follows

$$(3.8) \qquad \cdots \longrightarrow 0 \longrightarrow \pi_3(\Omega_0^\infty \mathrm{MTSO}(2)) \longrightarrow \mathbb{Z}/24$$
$$\hookrightarrow 0 \longrightarrow \pi_2(\Omega_0^\infty \mathrm{MTSO}(2)) \longrightarrow \pi_2^s(\mathbb{CP}^\infty) \oplus \mathbb{Z}/2 \xrightarrow[e']{d'} \mathbb{Z}/24.$$

There are natural maps from corresponding terms in 3.7 to that of 3.8. In Lemma 3.11 below, we will prove that the map $d$ has to be zero. Given Lemma 3.11, we obtain

$$(3.9) \qquad \pi_3(\Omega_0^\infty \mathrm{MT}\nu) \longrightarrow \pi_3(\Omega_0^\infty \mathrm{MTSO}(2)) = \mathbb{Z}/24.$$

The homology of $\Omega_0^\infty \mathrm{MTSO}(2)$ is the same as the stable homology of the mapping class group. Since the second stable homology of the mapping class group is $\mathbb{Z}$ ([Har83]), by the Hurewicz theorem we obtain that $\pi_2(\Omega_0^\infty \mathrm{MTSO}(2)) = \mathbb{Z}$. Hence, by comparing the maps $e$ and $e'$, we deduce that $\ker(e) = \mathbb{Z}$. Therefore, we have the short exact sequence

$$(3.10) \qquad 0 \to \pi_4^s(\overline{\mathrm{BS}\Gamma_2}) \to \pi_2(\Omega_0^\infty \mathrm{MT}\nu) \to \mathbb{Z} \to 0,$$

where the last map is induced by the map $\Omega_0^\infty \mathrm{MT}\nu \to \Omega_0^\infty \mathrm{MTSO}(2)$ on the second homotopy groups. Since $\overline{\mathrm{BS}\Gamma_2}$ is 3-connected, by the Hurewicz theorem, the fourth stable homotopy group of $\overline{\mathrm{BS}\Gamma_2}$ is the same as its fourth homotopy group, hence we obtain the second equality of the lemma. □

**Lemma 3.11.** *The map $d$ in the long exact sequence 3.7 is zero.*

*Proof.* Let $Q$ denote the functor $\Omega^\infty \Sigma^\infty$. To the circle bundle

$$\overline{\mathrm{BS}\Gamma_2} \to \mathrm{BS}\Gamma_2,$$

one can associate the circle transfer (i.e. the pretransfer for a circle bundle; see Remark 2.10) which is a map

$$\tau : Q((\mathrm{BS}\Gamma_2)_+) \to QS^{-1}((\overline{\mathrm{BS}\Gamma_2})_+),$$

where $QS^{-1}$ denotes the functor $\Omega^\infty \Sigma^{\infty-1}$. Recall the fiber sequence 3.4

$$(3.12) \qquad \Omega^\infty \mathrm{MT}\nu \to Q((\mathrm{BS}\Gamma_2)_+) \xrightarrow{\tau} QS^{-1}(\overline{\mathrm{BS}\Gamma_{2+}})$$

where the second map is the circle transfer $\tau$ by the construction of the fiber sequence. Hence the map $d$ is the map induced by the circle transfer $\tau$ on the third homotopy groups. In order to show that $d$ is zero, we consider the following pullback diagram

$$
\begin{array}{ccc}
S^1 & \longrightarrow & \overline{\mathrm{BS}\Gamma_2} \\
\downarrow & & \downarrow \\
* & \longrightarrow & \mathrm{BS}\Gamma_2,
\end{array}
$$

where the bottom horizontal map is a base point of $\mathrm{BS}\Gamma_2$. By naturality of the circle transfer, we have the following commutative diagram

$$
\begin{array}{ccc}
QS^{-1}(S_+^1) & \xrightarrow{g} & QS^{-1}(\overline{\mathrm{BS}\Gamma_{2+}}) \\
\uparrow {\scriptstyle f} & & \uparrow {\scriptstyle \tau} \\
QS^0 & \longrightarrow & Q((\mathrm{BS}\Gamma_2)_+),
\end{array}
$$

where $f$ is the circle transfer for the trivial circle bundle over a point and $g$ is the induced by naturality of the transfer map. Given that the map $\nu : \mathrm{BS}\Gamma_2 \to \mathrm{BGL}_2^+(\mathbb{R}) \simeq \mathbb{C}P^\infty$ is 4-connected (Remark 1.5), we obtain that

$$\pi_3^s(\mathrm{BS}\Gamma_2) = \pi_3^s(\mathbb{C}P^\infty) = 0.$$

Therefore, the bottom map in the above diagram induces an isomorphism on $\pi_3$. Hence to show that $d$ is zero, it is enough to show that $g$ induces zero on $\pi_3$. Note that the disjoint base point splits off naturally i.e. $QS^{-1}(S_+^1) \simeq QS^{-1}(S^1) \times QS^{-1}$ and similarly we have $QS^{-1}(\overline{\mathrm{B}S\Gamma_{2+}}) \simeq QS^{-1}(\overline{\mathrm{B}S\Gamma_2}) \times QS^{-1}$. Recall that $\pi_3(QS^{-1}) = \pi_4^s = 0$, so we need to show that the induced map

$$\pi_3(QS^{-1}(S^1)) \to \pi_3(QS^{-1}(\overline{\mathrm{B}S\Gamma_2}))$$

is zero. Given that $\overline{\mathrm{B}S\Gamma_2}$ is 3-connected (Remark 1.5) and in particular simply connected, the map $S^1 \to \overline{\mathrm{B}S\Gamma_2}$ is null-homotopic. Therefore, the map

$$QS^{-1}(S^1) \to QS^{-1}(\overline{\mathrm{B}S\Gamma_2}),$$

is also nullhomotopic which implies that the map $g$ induces zero on $\pi_3$.    □

**Theorem 3.13.** *Fix an embedding of $\mathbb{R}^2 \hookrightarrow \Sigma_{g,k}$. This embedding induces a map*

$$\mathrm{BDiff}_c^\delta(\mathbb{R}^2) \to \mathrm{BDiff}^\delta(\Sigma_{g,k}, \partial)$$

*which for $g \geq 4$ gives the following short exact sequence*

$$0 \to H_2(\mathrm{BDiff}_c^\delta(\mathbb{R}^2); \mathbb{Z}) \to H_2(\mathrm{BDiff}^\delta(\Sigma_{g,k}, \partial); \mathbb{Z}) \to H_2(\mathrm{BDiff}(\Sigma_{g,k}, \partial); \mathbb{Z}) \to 0.$$

*Remark* 3.14. Using Theorem 1.3 for $M \cong \mathbb{R}^2$ and the connectivity of $\overline{\mathrm{B}S\Gamma_2}$, one can in fact show that there is a natural map

$$H_2(\mathrm{BDiff}^\delta(\Sigma_{g,k}, \partial); \mathbb{Z}_{(p)}) \to H_4(\mathrm{B}S\Gamma_2; \mathbb{Z}_{(p)})$$

which induces isomorphism for $g \geq 4$ and the prime $p > 3$ (see [Nar]).

*Proof.* By Theorem 2.6, we know that for $g \geq 4$

$$H_2(\mathrm{BDiff}^\delta(\Sigma_{g,k}, \partial); \mathbb{Z}) \xrightarrow{\cong} H_2(\Omega_0^\infty \mathrm{MT}\nu; \mathbb{Z}).$$

Recall by Thurston's theorem on the perfectness of the identity component of diffeomorphism groups and Powell's theorem on the perfectness of the mapping class group, the group $\mathrm{Diff}(\Sigma_{g,k}; \partial)$ is a perfect group for $g > 2$ which implies that $\Omega_0^\infty \mathrm{MT}\nu$ is simply connected. Hence, by the Hurewicz theorem and 3.10, we have the following short exact sequence

$$(3.15) \qquad 0 \to \pi_4(\overline{\mathrm{B}S\Gamma_2}) \to H_2(\Omega_0^\infty \mathrm{MT}\nu) \to H_2(\Omega_0^\infty \mathrm{MTSO}(2)) \to 0.$$

By Theorem 1.3 or Remark 1.5, we know that $\overline{\mathrm{BDiff}_c(\mathbb{R}^2)}$ is homology equivalent to $\Omega^2 \overline{\mathrm{B}S\Gamma_2}$. Since $\mathrm{Diff}_c(\mathbb{R}^2)$ is contractible ([Sma59]), we have $\overline{\mathrm{BDiff}_c(\mathbb{R}^2)} \simeq \mathrm{BDiff}_c^\delta(\mathbb{R}^2)$. Given that $\overline{\mathrm{B}S\Gamma_2}$ is 3-connected, by the Hurewicz theorem we have $\pi_4(\overline{\mathrm{B}S\Gamma_2}) = H_2(\Omega^2 \overline{\mathrm{B}S\Gamma_2}; \mathbb{Z})$. Therefore, the exact sequence 3.15 is the same as the exact sequence in the theorem.    □

*Remark* 3.16. It is easy to see that the Serre spectral sequence for the fibration

$$\overline{\mathrm{BDiff}(\Sigma_g)} \to \mathrm{BDiff}^\delta(\Sigma_g) \to \mathrm{BDiff}(\Sigma_g)$$

and the perfectness of the identity component, $\mathrm{Diff}_0^\delta(\Sigma_g)$, implies that the map

$$H_2(\mathrm{BDiff}^\delta(\Sigma_g)) \longrightarrow H_2(\mathrm{BDiff}(\Sigma_g))$$

is surjective for all $g$, which means every surface bundle over a surface is cobordant to a flat surface bundle (see [KM05] for more explicit construction of such cobordisms).

We derive a geometric consequence of 3.10 and 3.9 for flat surface bundles.

**Theorem 3.17.** *Let $\mathrm{MSO}_k(X)$ denote the oriented cobordism group of $k$ dimensional manifolds equipped with a map to $X$. Then the map*

$$\mathrm{MSO}_3(\mathrm{BDiff}^\delta(\Sigma_g)) \longrightarrow\!\!\!\!\!\rightarrow \mathrm{MSO}_3(\mathrm{BDiff}(\Sigma_g))$$

*is surjective for $g \geq 6$. In other words, every surface bundle of genus at least $6$ over a 3-manifold is cobordant to a flat surface bundle.*

*Proof.* Note that $\mathrm{MSO}_i(X) = H_i(X; \mathbb{Z})$ for $i = 3$. Thus, we only need to prove it for homology. In order to prove that the following map is surjective in the stable range

$$H_3(\mathrm{BDiff}^\delta(\Sigma_{g,k}, \partial); \mathbb{Z}) \to H_3(\mathrm{BDiff}(\Sigma_{g,k}, \partial); \mathbb{Z})$$

we need a little lemma,

**Lemma 3.18.** *Let $X$ be a simply connected space, then we have the following exact sequence*

$$\pi_3(X) \to H_3(X) \to H_3(\mathrm{K}(\pi_2(X), 2); \mathbb{Z})$$

*where $\mathrm{K}(\pi_2(X), 2)$ is the Eilenberg-MacLane space whose second homotopy group is $\pi_2(X)$.*

The proof of the lemma is an easy Serre spectral sequence argument for the following map turned into fibration

$$X \to \mathrm{K}(\pi_2(X), 2).$$

Recall that since $\pi_2(\Omega_0^\infty \mathrm{MTSO}(2)) = \mathbb{Z}$ (see [Har83]), we can deduce that $H_3(\mathrm{K}(\pi_2(\Omega_0^\infty \mathrm{MTSO}(2)), 2); \mathbb{Z}) = 0$. If we apply lemma 3.18 for $X = \Omega_0^\infty \mathrm{MTSO}(2)$, we obtain that every degree 3 homology class of $\Omega_0^\infty \mathrm{MTSO}(2)$ is spherical i.e.

$$\pi_3(\Omega_0^\infty \mathrm{MTSO}(2)) \longrightarrow H_3(\Omega_0^\infty \mathrm{MTSO}(2); \mathbb{Z}).$$

We have the following commutative diagram by naturality of Hurewicz maps

$$
\begin{array}{ccc}
\pi_3(\Omega_0^\infty \mathrm{MT}\nu) & \longrightarrow & H_3(\Omega_0^\infty \mathrm{MT}\nu; \mathbb{Z}) \\
\downarrow & & \downarrow \\
\pi_3(\Omega_0^\infty \mathrm{MTSO}(2)) & \longrightarrow & H_3(\Omega_0^\infty \mathrm{MTSO}(2); \mathbb{Z}).
\end{array}
$$

the left vertical map is surjective by comparing exact sequences of 3.8 and 3.7, so the right vertical map has to be surjective. $\qquad \square$

### 3.2. Applications to characteristic classes of flat surface bundles.

In this section, we study two types of characteristic classes for flat surface bundles. The first type is constructed by forgetting the flat structure on the total space of the bundle and consider it just as a surface bundle. The second type of classes is secondary characteristic classes or so called Godbillon-Vey classes of the codimension 2 foliation induced by the flat structure on the total space.

#### 3.2.1. *The MMM-classes of the flat surface bundles.* Consider the following natural map again

$$\iota : \mathrm{BDiff}^\delta(\Sigma_g) \to \mathrm{BDiff}(\Sigma_g).$$

The first type of characteristic classes of the flat surface bundles is the pull back of MMM-classes via the map $\iota$. Recall the definition of MMM-classes is as follows. Let $\pi : E \to B$ be a surface bundle whose fibers are diffeomorphic to $\Sigma_g$. Let $\mathrm{T}\pi$ denote the vertical tangent bundle and let $e(\mathrm{T}\pi)$ denote its Euler class. Then, the i-th MMM class is defined to be

$$\kappa_i(E \to B) := \pi_!(e(\mathrm{T}\pi)^{i+1}) \in H^{2i}(B; \mathbb{Z})$$

where $\pi_!$ is the push forward map or the integration along the fiber which is defined since the fibers are compact closed manifolds. We denote the the i-th MMM class of the universal surface bundle over $\mathrm{BDiff}(\Sigma_g)$ by $\kappa_i$. Let $\kappa_i^\delta$ denote the pull back of $\kappa_i$ via the map $\iota$. One of the consequences of the Madsen-Weiss theorem and Harer stability is that the natural map

$$\mathbb{Q}[\kappa_1, \kappa_2, \dots] \to H^*(\mathrm{BDiff}(\Sigma_g); \mathbb{Q}),$$

is injective in the stable range. However, Morita in [Mor87] observed that $\kappa_i^\delta$ vanishes in rational cohomology for $i > 2$, i.e. the map

$$\mathbb{Q}[\kappa_1^\delta, \kappa_2^\delta, \dots] \to H^*(\mathrm{BDiff}^\delta(\Sigma_g); \mathbb{Q}),$$

sends $\kappa_i^\delta$ to zero while $i > 2$. For the above observation, it is essential to work with diffeomorphisms that are at least two times differentiable. It follows from Tsuboi's theorem ([Tsu89]) mentioned in the introduction and the Madsen-Weiss theorem that the similar map for $C^1$-diffeomorphisms is in fact an isomorphism.

Morita and Kotschick in [KM05] proved that there exists a flat surface bundle over a surface with nonzero signature, hence they conclude that $\kappa_1^\delta$ does not vanish in $H^2(\mathrm{BDiff}^\delta(\Sigma_g); \mathbb{Q})$ for $g \geq 3$. First, let us give a homotopy theoretic proof of their theorem using Theorem 2.6,

**Proposition 3.19.** $\kappa_1^\delta$ *does not vanish in* $H^2(\mathrm{BDiff}^\delta(\Sigma_g); \mathbb{Q})$ *for* $g > 3$.

*Proof.* By Theorem 2.6 and Theorem 1.14, we know

$$H^2(\Omega_0^\infty \mathrm{MT}\nu; \mathbb{Q}) \to H^2(\mathrm{BDiff}^\delta(\Sigma_g); \mathbb{Q}),$$

is surjective for $g \geq 3$ and an isomorphism for $g > 3$. Thus, we need to show that the corresponding class in $H^2(\Omega_0^\infty \mathrm{MT}\nu; \mathbb{Q})$ is nonzero. Consider the following sequence

$$H^{2i+2}(\mathrm{BS}\Gamma_2; \mathbb{Q}) \xrightarrow{\cong} H^{2i}(\mathrm{MT}\nu; \mathbb{Q}) \xrightarrow{\sigma^*} H^{2i}(\Omega_0^\infty \mathrm{MT}\nu; \mathbb{Q}),$$

where the first map is the Thom isomorphism and the second map is induced by suspension map. Let $e$ denote the generator of $H^2(\mathrm{BSO}(2); \mathbb{Q})$. It is easy to see that the image of $\nu^*(e^{i+1}) \in H^{2i+2}(\mathrm{BS}\Gamma_2; \mathbb{Q})$ in $H^{2i}(\Omega_0^\infty \mathrm{MT}\nu; \mathbb{Q})$ is $\kappa_i^\delta$ (see [GMT06, Theorem 3.1]). Since the map $\nu$ is at least 4-connected, $\nu^*(e^2)$ is not zero. Thus, to prove that $\kappa_1^\delta$ is nonzero, it is enough to show that the following map is nontrivial

$$H_2(\Omega_0^\infty \mathrm{MT}\nu; \mathbb{Q}) \xrightarrow{\sigma_*} H_2(\mathrm{MT}\nu; \mathbb{Q}) \xrightarrow{\kappa_1^\delta} \mathbb{Q}.$$

To prove that the suspension map is surjective on rational homology, let us consider the following commutative diagram

$$\begin{array}{ccc}
\pi_*(\Omega^\infty \mathrm{MT}\nu) \otimes \mathbb{Q} & \longrightarrow & \pi_*(\mathrm{MT}\nu) \otimes \mathbb{Q} \\
\downarrow & & \downarrow \\
H_*(\Omega^\infty \mathrm{MT}\nu; \mathbb{Q}) & \longrightarrow & H_*(\mathrm{MT}\nu; \mathbb{Q}).
\end{array}$$

The horizontal maps are induced by the suspension map and the vertical maps are induced by the Hurewicz map. The top horizontal map is an isomorphism by the definition of the homotopy groups of a spectra and the right vertical map is also an isomorphism because of the rational Hurewicz theorem (see [Rud98, Theorem 7.11]). Therefore, $\sigma_*$, the bottom horizontal map, is surjective, which implies $\kappa_1^\delta \circ \sigma_*$ is nontrivial. □

*Remark* 3.20. Morita and Kotschick in [KM05] by a formal argument showed that nontriviality of $\kappa_1^\delta$ implies that all its powers are nontrivial in the stable range. This can also be deduced from the fact that $H_*(\Omega^\infty \mathrm{MT}\nu; \mathbb{Q})$ is a Hopf algebra over $\mathbb{Q}$.

Regarding $\kappa_2^\delta$, Morita and Kotschick asked the following problem

**Problem.** *Does* $\kappa_2^\delta$ *vanish in* $H^4(\mathrm{BDiff}^\delta(\Sigma_g); \mathbb{Q})$?

Toward answering this problem, we prove that it is equivalent to an open problem in foliation theory.

**Theorem 3.21.** *The MMM-class $\kappa_2$ in $H^4(\mathrm{BDiff}^\delta(\Sigma_g);\mathbb{Q})$ is nonzero for $g \geq 6$ if and only if a $C^2$-foliation $\mathcal{F}$ of codimension 2 on a 6 manifold exists such that $e(\nu(\mathcal{F}))^3 \neq 0$ where $\nu(\mathcal{F})$ is the normal bundle of the foliation $\mathcal{F}$.*

*Proof.* Let $\mathcal{F}\Omega_{n,k}$ denote the group of foliated cobordism group of codimension $k$ foliations on $n$ dimensional manifolds. Every oriented codimension $k$ foliation on a manifold $M$ gives a well-defined homotopy class of maps from $M$ to $\mathrm{B}\Gamma_k$ (see [Hae71, Theorem 7]). Hence, we have the well defined map

$$\mathcal{F}\Omega_{n,k} \to \mathrm{MSO}_n(\mathrm{B}\Gamma_k).$$

By a result of [Kos74, Theorem 1'], the map $\mathcal{F}\Omega_{6,2} \to \mathrm{MSO}_6(\mathrm{B}\Gamma_2)$ is rationally bijective. By Atiyah-Hirzebruch spectral sequence $\mathrm{MSO}_4(\mathrm{BDiff}^\delta(\Sigma_g)) \to H_4(\mathrm{BDiff}^\delta(\Sigma_g))$ is surjective. Hence, in order to prove (non)triviality of $\kappa_2^\delta$, we need to study (non)triviality of the map

$$\kappa_2^\delta \colon \mathrm{MSO}_4(\mathrm{BDiff}^\delta(\Sigma_g)) \to \mathbb{Q}.$$

By Theorem 2.6, we know

$$\mathrm{MSO}_4(\mathrm{BDiff}^\delta(\Sigma_g)) \to \mathrm{MSO}_4(\Omega_0^\infty \mathrm{MT}\nu)$$

is an isomorphism for $g > 6$ and surjective for $g \geq 6$. We have the following commutative diagram

$$
\begin{array}{ccc}
\mathrm{MSO}_4(\Omega_0^\infty \mathrm{MT}\nu) \otimes \mathbb{Q} & \xrightarrow{\ \kappa_2^\delta\ } & \mathbb{Q} \\
\downarrow{\scriptstyle \sigma_*} & & \uparrow{\scriptstyle e^3} \\
\mathrm{MSO}_4(\mathrm{MT}\nu) \otimes \mathbb{Q} & \xrightarrow{\ \cong\ } & \mathrm{MSO}_6(\mathrm{B}\Gamma_2) \otimes \mathbb{Q}
\end{array}
$$

where $\sigma$ is the suspension map, the bottom map is the Thom isomorphism and the right vertical map is the map given by the cube of the Euler class of the codimension 2 Haefliger structures. Hence, if we show that the map $\sigma_*$ is rationally surjective, (non)triviality of $\kappa_2^\delta$ and $e^3$ become equivalent. We know that

$$H_4(\Omega_0^\infty \mathrm{MT}\nu;\mathbb{Q}) \to H_4(\mathrm{MT}\nu;\mathbb{Q})$$

is surjective by the same argument as in proposition 3.19 and since the Atiyah-Hirzebruch spectral sequence implies that $\mathrm{MSO}_4(X) = H_4(X) \oplus \mathbb{Z}$, we have

$$\sigma_* \colon \mathrm{MSO}_4(\Omega_0^\infty \mathrm{MT}\nu) \otimes \mathbb{Q} \twoheadrightarrow \mathrm{MSO}_4(\mathrm{MT}\nu) \otimes \mathbb{Q}$$

$\square$

*Remark* 3.22. One possible way to show that $e^3$ is nonzero in $H^6(\mathrm{B}\Gamma_2;\mathbb{Q})$ which implies that a flat surface bundle with a nontrivial $\kappa_2^\delta$ exists, is to look at the Gysin sequence for the circle bundle

$$S^1 \to \overline{\mathrm{B}\Gamma_2} \to \mathrm{B}\Gamma_2.$$

Part of the Gysin sequence that is relevant for us is

$$H_6(\mathrm{B}\Gamma_2;\mathbb{Q}) \xrightarrow{\ \cap e\ } H_4(\mathrm{B}\Gamma_2;\mathbb{Q}) \xrightarrow{\ \tau\ } H_5(\overline{\mathrm{B}\Gamma_2};\mathbb{Q}),$$

where $\tau$ is the transgression map for the circle bundle. If we show that a class $c \in H_4(\mathrm{B}\Gamma_2;\mathbb{Q})$ that satisfies $e^2(c) \neq 0$, maps to zero via the transgression $\tau$, then the Gysin sequence implies that $e^3$ has to be nonzero in $H^6(\mathrm{B}\Gamma_2;\mathbb{Q})$. To construct a class $c$ for which $e^2(c) \neq 0$, we consider a map $f\colon \mathbb{C}P^2 \to \mathrm{BGL}_2^+(\mathbb{R})$ that is nontrivial on the fourth homology with rational coefficients. Since $\mathrm{B}\Gamma_2 \to \mathrm{BGL}_2^+(\mathbb{R})$ is 4-connected (Remark 1.5), we can lift $f$ to a map $g\colon \mathbb{C}P^2 \to \mathrm{B}\Gamma_2$. If we take the pullback of the circle bundle $\overline{\mathrm{BS\Gamma}}_2 \to \mathrm{BS\Gamma}_2$ via the map $g$, we obtain the following commutative diagram

$$
\begin{array}{ccc}
S^5 & \xrightarrow{\hat{g}} & \overline{\mathrm{BS\Gamma}}_2 \\
\downarrow & & \downarrow \\
\mathbb{C}P^2 & \xrightarrow{g} & \mathrm{BS\Gamma}_2.
\end{array}
$$

The transgression maps the class $g([\mathbb{C}P^2]) \in H_4(\mathrm{BS\Gamma}_2; \mathbb{Q})$ to the class $\hat{g}([S^5]) \in H_5(\overline{\mathrm{BS\Gamma}}_2; \mathbb{Q}) \cong \pi_5(\overline{\mathrm{BS\Gamma}}_2) \otimes \mathbb{Q}$. Since $\pi_5(\overline{\mathrm{BS\Gamma}}_2) = \pi_5(\mathrm{BS\Gamma}_2)$, we conclude that $e^3$ is nontrivial in $H^6((\mathrm{BS\Gamma}_2; \mathbb{Q})$ if

$$
g_* : \pi_5(\mathbb{C}P^2) \otimes \mathbb{Q} \to \pi_5(\mathrm{BS\Gamma}_2) \otimes \mathbb{Q}
$$

is trivial. Note that since the image of $g$ maps to zero in $H_5(\mathrm{BS\Gamma}_2; \mathbb{Q})$, the Godbillon-Vey classes vanish on the image of $g$, so it cannot be detected by the GV classes.

### 3.3. **MMM-classes as integral cohomology classes.**

The situation, however, is surprisingly different with integer coefficients. To study integral MMM-classes, we first reduce them to classes with finite field coefficients. Consider the following commutative diagram

$$
\begin{array}{ccc}
H^*(\Omega_0^\infty \mathrm{MTSO}(2); \mathbb{F}_p) & \xrightarrow{\Omega^\infty \boldsymbol{\nu}^*} & H^*(\Omega_0^\infty \mathrm{MT}\nu; \mathbb{F}_p) \\
\downarrow & & \downarrow \\
H^*(\mathrm{BDiff}(\Sigma_g); \mathbb{F}_p) & \xrightarrow{\iota^*} & H^*(\mathrm{BDiff}^\delta(\Sigma_g); \mathbb{F}_p).
\end{array}
$$

where the vertical maps are isomorphisms in the stable range. Thus, to study $\iota^*$, we need to study $\Omega^\infty \boldsymbol{\nu}^*$. In Theorem 2.12, we proved that $\Omega^\infty \boldsymbol{\nu}^*$ is injective on cohomology with $\mathbb{F}_p$ coefficients. Given the injectivity on $\mathbb{F}_p$-cohomology for all prime $p$, as we shall see the map

$$
H^*(\Omega_0^\infty \mathrm{MTSO}(2); \mathbb{Z}) \longrightarrow H^*(\Omega_0^\infty \mathrm{MT}\nu; \mathbb{Z})
$$

is injective. Hence, we summarize the situation with finite coefficients and with integer coefficients as follows

**Theorem 3.23.** *For every prime p, the map*

$$
H^*(\mathrm{BDiff}(\Sigma_g); \mathbb{F}_p) \xrightarrow{\iota^*} H^*(\mathrm{BDiff}^\delta(\Sigma_g); \mathbb{F}_p),
$$

*is injective in the stable range.*

**Theorem 3.24.** *The induced map*

$$
\mathbb{Z}[\kappa_1^\delta, \kappa_2^\delta, \dots] \longhookrightarrow H^*(\mathrm{BDiff}^\delta(\Sigma_{\infty,1}, \partial); \mathbb{Z}),
$$

*is injective.*

*Proof.* We prove even a stronger result that the induced map

$$
\iota^* : H^*(\mathrm{BDiff}(\Sigma_{\infty,1}, \partial); \mathbb{Z}) \to H^*(\mathrm{BDiff}^\delta(\Sigma_{\infty,1}, \partial); \mathbb{Z}),
$$

is injective. It follows from Theorem 2.12 that the map

$$
\iota_p^* : H^*(\mathrm{BDiff}(\Sigma_{\infty,1}, \partial); \mathbb{F}_p) \to H^*(\mathrm{BDiff}^\delta(\Sigma_{\infty,1}, \partial); \mathbb{F}_p),
$$

is injective for all prime $p$. To prove that $\iota^*$ is injective, consider the induced map between the Bockstein exact sequences

$$H^*(\mathrm{BDiff}(\Sigma_{\infty,1},\partial);\mathbb{F}_p) \longleftarrow H^*(\mathrm{BDiff}(\Sigma_{\infty,1},\partial);\mathbb{Z}) \overset{\times p}{\longleftarrow} H^*(\mathrm{BDiff}(\Sigma_{\infty,1},\partial);\mathbb{Z})$$

$$\downarrow \iota_p^* \qquad\qquad\qquad \downarrow \iota^* \qquad\qquad\qquad \downarrow \iota^*$$

$$H^*(\mathrm{BDiff}^\delta(\Sigma_{\infty,1},\partial);\mathbb{F}_p) \longleftarrow H^*(\mathrm{BDiff}^\delta(\Sigma_{\infty,1},\partial);\mathbb{Z}) \overset{\times p}{\longleftarrow} H^*(\mathrm{BDiff}^\delta(\Sigma_{\infty,1},\partial);\mathbb{Z}).$$

Let $a \in \mathrm{Ker}(\iota^*)$, since $\iota_p^*$ is injective for all $p$, it follows from the above diagram that $a \in pH^*(\mathrm{BDiff}(\Sigma_{\infty,1},\partial);\mathbb{Z})$ for all $p$. Since $H^*(\mathrm{BDiff}(\Sigma_{\infty,1},\partial);\mathbb{Z})$ is finitely generated by the Madsen-Weiss theorem ([MW07]), we deduce $a = 0$. □

*Remark* 3.25. Akita, Kawazumi and Uemura in [AKU01] proved the algebraic independence of MMM-classes by using finite cyclic subgroups of the mapping class groups. Their method can be also applied to prove the above theorem.

*Remark* 3.26. Recall that as Morita in [Mor87] observed, the classes $\kappa_i^\delta$ for $i > 2$ in cohomology with $\mathbb{R}$ coefficients or even with $\mathbb{Q}$ coefficients vanish. This observation implies that there is a class in $H^{2i-1}(\mathrm{BDiff}^\delta(\Sigma_g);\mathbb{R}/\mathbb{Z})$ for $i > 2$ which maps to $\kappa_i^\delta$ in the Bockstein exact sequence

$$H^{2i-1}(\mathrm{BDiff}^\delta(\Sigma_g);\mathbb{R}/\mathbb{Z}) \to H^{2i}(\mathrm{BDiff}^\delta(\Sigma_g);\mathbb{Z}) \to H^{2i}(\mathrm{BDiff}^\delta(\Sigma_g);\mathbb{R}).$$

Cheeger-Simons character theory helps us to find a canonical lift of $\kappa_{2i+1}^\delta$ for $i > 0$ in $H^{4i+1}(\mathrm{BDiff}^\delta(\Sigma_g);\mathbb{R}/\mathbb{Z})$. In [CS85, Proposition 7.3], Cheeger and Simons showed that there are canonical classes known as Pontryagin characters $\widehat{p_i} \in H^{4i-1}(\mathrm{B\Gamma}_n;\mathbb{R}/\mathbb{Z})$ for $2i > n$ so that the image of $\widehat{p_i}$ under Bockstein map

$$H^{4i-1}(\mathrm{B\Gamma}_n;\mathbb{R}/\mathbb{Z}) \overset{\beta}{\longrightarrow} H^{4i}(\mathrm{B\Gamma}_n;\mathbb{Z})$$

is $-p_i$ of the normal bundle of the universal Haefliger structure. Let $\widehat{\kappa_{2i+1}}$ denote the image of the class $\widehat{p_1^{i+1}} \in H^{4i+3}(\mathrm{B\Gamma}_2;\mathbb{R}/\mathbb{Z})$ under the following maps

$$H^{4i+3}(\mathrm{B\Gamma}_2;\mathbb{R}/\mathbb{Z}) \to H^{4i+1}(\mathrm{MT}\nu;\mathbb{R}/\mathbb{Z}) \to H^{4i+1}(\Omega_0^\infty \mathrm{MT}\nu;\mathbb{R}/\mathbb{Z}).$$

Hence, the class $\widehat{\kappa_{2i+1}}$ gives a canonical lift of $\kappa_{2i+1}^\delta$ if their degree lie in the stable range and $\beta(\widehat{\kappa_{2i+1}}) = -\kappa_{2i+1}^\delta$ where $\beta$ is the Bockstein map in the following sequence

$$H^{4i+1}(\mathrm{BDiff}^\delta(\Sigma_g);\mathbb{R}/\mathbb{Z}) \overset{\beta}{\longrightarrow} H^{4i+2}(\mathrm{BDiff}^\delta(\Sigma_g);\mathbb{Z}) .$$

We showed that $\kappa_{2i+1}^\delta$ is a non-torsion class so $\beta$ is a nontrivial map. Therefore, for $i > 0$ we obtain a regulator type map for discrete surface diffeomorphisms

$$\text{regulator map} : H_{4i+1}(\mathrm{BDiff}^\delta(\Sigma_g);\mathbb{Z}) \overset{\widehat{\kappa_{2i+1}}}{\longrightarrow} \mathbb{R}/\mathbb{Z} .$$

*Question.* What is the cocycle formula for $\widehat{\kappa_{2i+1}}$?

Nontriviality of this regulator map implies that $H_{4i+1}(\mathrm{BDiff}^\delta(\Sigma_g);\mathbb{Z})$ in the stable range is not a trivial group but in fact it follows easily from Theorem 3.24 that $H_{2i+1}(\mathrm{BDiff}^\delta(\Sigma_g);\mathbb{Z})$ in the stable range and for $i > 1$ is in fact uncountable (see [Nar14, Theorem 6.4]).

3.3.1. *Secondary classes of flat surface bundles.* In this section, we prove Theorem 0.18 in the introduction. Let us first recall a preliminary result about continuous variation of secondary characteristic classes of foliations of codimension 2. The cohomology of $\mathrm{B\Gamma}_n$ is not yet very well understood but it has been extensively studied via secondary characteristic classes of foliations known as the Godbillon-Vey classes. For codimension 2 foliations there are two GV-classes denoted by $h_1c_2, h_1c_1^2 \in H^5(\mathrm{B\Gamma}_2; \mathbb{R})$ (for definition of these classes look at [Pit76], [Law77]). Rasmussen in [Ras80] proved that these two classes vary continuously and independently i.e.

$$(h_1c_2, h_1c_1^2) : H_5(\mathrm{B\Gamma}_2; \mathbb{Z}) \longrightarrow \mathbb{R}^2.$$

If we take the universal flat surface bundle over $\mathrm{BDiff}^\delta(\Sigma_g)$, we can integrate $h_1c_2$ and $h_1c_1^2$ along the compact fibers and we denote their integration along the fiber as

$$\int_{\text{fiber}} h_1c_2, \int_{\text{fiber}} h_1c_1^2 \in H^3(\mathrm{BDiff}^\delta(\Sigma_g); \mathbb{R}).$$

Morita in [Mor06, Problem 44] posed the question that whether the map induced by the above two cohomology classes from $H_3(\mathrm{BDiff}^\delta(\Sigma_g); \mathbb{Z})$ to $\mathbb{R}^2$ is surjective. Bowden in [Bow12] used a curious spectral sequence that only converges in low homological degrees to answer Morita's question affirmatively. Here, we simplify his proof using the Mather-Thurston theorem.

**Theorem 3.27** (Bowden). *For any $g$, the induced map*

$$k : H_3(\mathrm{BDiff}^\delta(\Sigma_g); \mathbb{Q}) \xrightarrow{\;(\int_{\text{fiber}} h_1c_2, \int_{\text{fiber}} h_1c_1^2)\;} \mathbb{R}^2$$

*is surjective.*

*Proof.* Embed $\mathbb{R}^2$ as the interior of a small disk into $\Sigma_g$. The restriction of the map $k$ to the embedded disk gives the following commutative diagram

$$
\begin{array}{ccc}
H_3(\mathrm{BDiff}^\delta(\Sigma_g); \mathbb{Q}) & \xrightarrow{\;\;k\;\;} & \mathbb{R}^2 \\
\uparrow & \nearrow{\scriptstyle k'} & \\
H_3(\mathrm{BDiff}_c^\delta(\mathbb{R}^2); \mathbb{Q}), & &
\end{array}
$$

where the map $k'$ is also induced by the fiber integration of the GV classes along the embedded disk. If we show $k'$ is surjective then it implies that $k$ is also surjective. Let $\overline{\mathrm{BDiff}_c(\mathbb{R}^2)}$ be the homotopy fiber of the map

$$\mathrm{BDiff}_c^\delta(\mathbb{R}^2) \to \mathrm{BDiff}_c(\mathbb{R}^2).$$

But the topological group $\mathrm{Diff}_c(\mathbb{R}^2)$ is contractible ([Sma59]), so $\overline{\mathrm{BDiff}_c(\mathbb{R}^2)} \simeq \mathrm{BDiff}_c^\delta(\mathbb{R}^2)$. By Thurston's theorem [Thu74], we know that there is a map

$$\overline{\mathrm{BDiff}_c(\mathbb{R}^2)} \to \Omega^2 \overline{\mathrm{B\Gamma}_2}$$

that induces a homology isomorphism. By Remark 1.5, we know $\overline{\mathrm{B\Gamma}_2}$ is at least 3-connected. Therefore, we have

$$(3.28) \qquad H_3(\mathrm{BDiff}_c^\delta(\mathbb{R}^2); \mathbb{Q}) \xrightarrow{\cong} H_3(\Omega^2 \overline{\mathrm{B\Gamma}_2}; \mathbb{Q}) \twoheadrightarrow H_5(\overline{\mathrm{B\Gamma}_2}; \mathbb{Q}) \twoheadrightarrow \mathbb{R}^2.$$

The first map is an isomorphism by Thurston's theorem. The second map is the suspension map and because $\overline{\mathrm{B\Gamma}_2}$ is at least 3-connected, the rational Hurewicz theorem implies that in the diagram

$$
\begin{array}{ccc}
\pi_3(\Omega^2\overline{\mathrm{BST}_2}) \otimes \mathbb{Q} & \longrightarrow & H_3(\Omega^2\overline{\mathrm{BST}_2}; \mathbb{Q}) \\
\downarrow{\scriptstyle\cong} & & \downarrow \\
\pi_5(\overline{\mathrm{BST}_2}) \otimes \mathbb{Q} & \longrightarrow & H_5(\overline{\mathrm{BST}_2}; \mathbb{Q}),
\end{array}
$$

the bottom map is surjective, so is the the left vertical map. The third map in 3.28 is given by the Godbillon-Vey classes

$$
H_5(\overline{\mathrm{BST}_2}; \mathbb{Q}) \xrightarrow{(\int h_1 c_2, \int h_1 c_1^2)} \mathbb{R}^2
$$

which is surjective as a corollary of the theorem of Rasmussen [Ras80]. □

**Corollary 3.29.** *There exists a surjective map*

$$
H_{3k}(\mathrm{BDiff}^\delta(\Sigma_g); \mathbb{Q}) \longrightarrow \wedge^k_{\mathbb{Q}} \mathbb{R}^2,
$$

*provided $k \leq (2g-2)/9$ where $\wedge^k_{\mathbb{Q}} \mathbb{R}^2$ is the $k$-th exterior power of $\mathbb{R}^2$ as the vector space over $\mathbb{Q}$.*

*Proof.* For $k \leq (2g-2)/9$, using Theorem 2.6, we know that

$$
H_{3k}(\mathrm{BDiff}^\delta(\Sigma_g); \mathbb{Q}) \simeq H_{3k}(\Omega_0^\infty \mathrm{MT}\nu; \mathbb{Q}).
$$

By Theorem 3.27, we obtain a surjective map

$$
H_3(\Omega_0^\infty \mathrm{MT}\nu; \mathbb{Q}) \longrightarrow \mathbb{R}^2.
$$

Note that $H_*(\Omega_0^\infty \mathrm{MT}\nu; \mathbb{Q})$ is a simply connected Hopf algebra over $\mathbb{Q}$. Hence elements in $H_3(\Omega_0^\infty \mathrm{MT}\nu; \mathbb{Q})$ are primitive in the Hopf algebra. If we choose a basis for the vector space $H_3(\Omega_0^\infty \mathrm{MT}\nu; \mathbb{Q})$, their exterior powers are nontrivial and provide us with a map

$$
H_{3k}(\Omega_0^\infty \mathrm{MT}\nu; \mathbb{Q}) \longrightarrow \wedge^k_{\mathbb{Q}} \mathbb{R}^2,
$$

which is surjective. □

Recall from [H+91] that $H_3(\mathrm{BDiff}(\Sigma_g); \mathbb{Q}) = 0$, but as we showed in the proof of Theorem 3.27, secondary characteristic classes in $H_3(\mathrm{BDiff}^\delta(\Sigma_g); \mathbb{Q})$ vary continuously and independently on diffeomorphisms of $\Sigma_g$ that are only supported in a disk. Bowden asked the author if there is a nontrivial class in $H_3(\mathrm{BDiff}^\delta(\Sigma_g); \mathbb{Q})$ that cannot be detected by an embedding of a disk. To give an answer to his question, we prove that at least in the stable range, all classes in $H_3(\mathrm{BDiff}^\delta(\Sigma_g); \mathbb{Q})$ are essentially supported in a disk, more precisely

**Theorem 3.30.** *Let $\mathbb{R}^2 \hookrightarrow \Sigma_g$ be an embedding of an open disk in the surface $\Sigma_g$. For $g \geq 6$, the induced map*

$$
H_3(\mathrm{BDiff}_c^\delta(\mathbb{R}^2); \mathbb{Q}) \longrightarrow H_3(\mathrm{BDiff}^\delta(\Sigma_g); \mathbb{Q}),
$$

*is surjective.*

*Proof.* Recall that since the topological group $\mathrm{Diff}_c(\mathbb{R}^2)$ is contractible (see [Sma59]), we have $\overline{\mathrm{BDiff}_c(\mathbb{R}^2)} \simeq \mathrm{BDiff}_c^\delta(\mathbb{R}^2)$. Hence, all the $\mathbb{R}^2$-bundles trivialized at the infinity over $\mathrm{BDiff}_c^\delta(\mathbb{R}^2)$ are topologically trivial bundle, therefore Pontryagin-Thom theory for the trivial bundle $\mathrm{BDiff}_c^\delta(\mathbb{R}^2) \times \mathbb{R}^2 \to \mathrm{BDiff}_c^\delta(\mathbb{R}^2)$, as we shall explain in Lemma 3.35, provides a map

$$
\beta: \mathrm{BDiff}_c^\delta(\mathbb{R}^2) \to \Omega^2 \overline{\mathrm{BST}_2},
$$

which is a homology isomorphism by the theorem of Thurston ([Thu74]). We showed that there exists a Madsen-Weiss type map

$$\alpha_{\Sigma_g} : \mathrm{BDiff}^\delta(\Sigma_g) \to \Omega_0^\infty \mathrm{MT}\nu,$$

for the universal flat $\Sigma_g$-bundle over $\mathrm{BDiff}^\delta(\Sigma_g)$ that induces a homology isomorphism in the stable range. In Lemma 3.35 below, we shall prove that there exists a homotopy commutative diagram

$$(3.31) \quad \begin{array}{ccc}
& \mathrm{BDiff}_c^\delta(\mathbb{R}^2) \longrightarrow \mathrm{BDiff}^\delta(\Sigma_g) \\
\beta \swarrow & \downarrow \alpha_{\mathbb{R}^2} & \downarrow \alpha_{\Sigma_g} \\
\Omega^2 \overline{\mathrm{B}S\Gamma_2} \xrightarrow{\lambda} Q_0 S^{-2}(\overline{\mathrm{B}S\Gamma_{2+}}) \longrightarrow \Omega_0^\infty \mathrm{MT}\nu,
\end{array}$$

where $Q = \Omega^\infty \Sigma^\infty$ and the subscript $0$ means the base point component. The maps $\alpha_{\mathbb{R}^2}$ and $\alpha_{\Sigma_g}$ are defined as Remark 2.10 and the map $\lambda$ is the natural composition

$$\Omega^2 \overline{\mathrm{B}S\Gamma_2} \to Q(\Omega^2 \overline{\mathrm{B}S\Gamma_2}) \to QS^{-2}(\overline{\mathrm{B}S\Gamma_{2+}}).$$

Thus to prove the theorem, we need to show that the induced map

$$H_3(\Omega^2 \overline{\mathrm{B}S\Gamma_2}; \mathbb{Q}) \to H_3(\Omega_0^\infty \mathrm{MT}\nu; \mathbb{Q})$$

is surjective. Consider the following commutative diagram

$$(3.32) \quad \begin{array}{ccc}
\pi_3(\Omega^2 \overline{\mathrm{B}S\Gamma_2}) \otimes \mathbb{Q} & \xrightarrow{g} & \pi_3(\Omega_0^\infty \mathrm{MT}\nu) \otimes \mathbb{Q} \\
\downarrow & & \downarrow \\
H_3(\Omega^2 \overline{\mathrm{B}S\Gamma_2}; \mathbb{Q}) & \xrightarrow{f} & H_3(\Omega_0^\infty \mathrm{MT}\nu; \mathbb{Q}).
\end{array}$$

Because we have $H_1(\Omega_0^\infty \mathrm{MT}\nu; \mathbb{Q}) = H_3(\mathrm{B}S\Gamma_2; \mathbb{Q}) = 0$, hence the rational Hurewicz theorem implies that the right vertical map is surjective. Therefore to prove that the map $f$ in 3.32 is surjective, it is sufficient to prove that the map $g$ is surjective.

Recall from Theorem 3.3 that we have the fibration sequence

$$QS^{-2}(\overline{\mathrm{B}S\Gamma_{2+}}) \xrightarrow{h} \Omega^\infty \mathrm{MT}\nu \to Q(\mathrm{B}S\Gamma_{2+}),$$

and its long exact sequence of homotopy groups in 3.7 implies that the map $h$ induces a surjection between third rational homotopy groups

$$(3.33) \quad \pi_5^s(\overline{\mathrm{B}S\Gamma_2}) \otimes \mathbb{Q} \xrightarrow{h_*} \pi_3(\Omega_0^\infty \mathrm{MT}\nu) \otimes \mathbb{Q}.$$

As for the surjectivity of the map $g$, we observed in the diagram 3.31 that the map $g$ is induced by the composition

$$(3.34)$$
$$\pi_3(\Omega^2 \overline{\mathrm{B}S\Gamma_2}) \otimes \mathbb{Q} = \pi_5(\overline{\mathrm{B}S\Gamma_2}) \otimes \mathbb{Q} \xrightarrow{\lambda_*} \pi_3(QS^{-2}(\overline{\mathrm{B}S\Gamma_{2+}})) \otimes \mathbb{Q} \xrightarrow{h_*} \pi_3(\Omega_0^\infty \mathrm{MT}\nu) \otimes \mathbb{Q}.$$

Note that $\pi_3(QS^{-2}(\overline{\mathrm{B}S\Gamma_{2+}})) = \pi_5^s(\overline{\mathrm{B}S\Gamma_2}) \oplus \pi_5^s$ and since $\pi_5^s = 0$, the map $g$ is surjective if in 3.34 the map $\lambda_*$ is surjective. Recall $\overline{\mathrm{B}S\Gamma_2}$ is 3-connected (Remark 1.5), therefore the rational Hurewicz theorem implies that the Hurewicz map

$$\pi_5(\overline{\mathrm{B}S\Gamma_2}) \otimes \mathbb{Q} \xrightarrow{\cong} \pi_5^s(\overline{\mathrm{B}S\Gamma_2}) \otimes \mathbb{Q},$$

is an isomorphism. Hence $g = h_* \circ \lambda_*$ is surjective.                    $\square$

**Lemma 3.35.** *The diagram 3.31 is homotopy commutative.*

*Proof.* Since the group $\mathrm{Diff}_c^\delta(\mathbb{R}^2)$ acts on the surface $\Sigma_g$ via the embedding $\mathbb{R}^2 \hookrightarrow \Sigma_g$, we obtain a map between Borel constructions

$$\mathbb{R}^2 /\!\!/ \mathrm{Diff}_c^\delta(\mathbb{R}^2) \to \Sigma_g /\!\!/ \mathrm{Diff}_c^\delta(\mathbb{R}^2).$$

On the other hand $\Sigma_g /\!\!/ \mathrm{Diff}_c^\delta(\mathbb{R}^2)$ is a flat surface bundle induced by the pullback of the universal flat surface bundle over $\mathrm{BDiff}^\delta(\Sigma_g)$ via the map $\mathrm{BDiff}_c^\delta(\mathbb{R}^2) \to \mathrm{BDiff}^\delta(\Sigma_g)$. Therefore, we have the homotopy commutative diagram

$$\begin{array}{ccccc}
\mathbb{R}^2 /\!\!/ \mathrm{Diff}_c^\delta(\mathbb{R}^2) & \longrightarrow & \Sigma_g /\!\!/ \mathrm{Diff}_c^\delta(\mathbb{R}^2) & \longrightarrow & \Sigma_g /\!\!/ \mathrm{Diff}^\delta(\Sigma_g) \\
& \searrow_{\pi'} & \downarrow & & \pi \downarrow \\
& & \mathrm{BDiff}_c^\delta(\mathbb{R}^2) & \longrightarrow & \mathrm{BDiff}^\delta(\Sigma_g).
\end{array}$$

By the naturality of the pretransfer construction in Remark 2.10, we have a commutative diagram of spectra

$$(3.36) \quad \begin{array}{ccc}
\Sigma^\infty(\mathrm{BDiff}_c^\delta(\mathbb{R}^2)_+) & \longrightarrow & \mathbf{Th}(-T\pi') \\
\downarrow & & \downarrow \\
\Sigma^\infty(\mathrm{BDiff}^\delta(\Sigma_g)_+) & \longrightarrow & \mathbf{Th}(-T\pi).
\end{array}$$

Note that by the flatness of the bundles the classifying map for $T\pi$ and $T\pi'$ lift to $\mathrm{B}S\Gamma_2$. But since $\mathrm{BDiff}_c(\mathbb{R}^2)$ is contractible, topologically the bundle $\mathbb{R}^2 /\!\!/ \mathrm{Diff}_c^\delta(\mathbb{R}^2)$ is trivial which implies that the bundle $T\pi'$ is trivial, hence the classifying map for $T\pi'$ further lifts to $\overline{\mathrm{B}S\Gamma_2}$. Thus, we have a commutative diagram

$$(3.37) \quad \begin{array}{ccc}
\mathbf{Th}(-T\pi') & \longrightarrow & \mathbf{Th}(-\underline{\mathbb{R}^2} \to \overline{\mathrm{B}S\Gamma_2}) \\
\downarrow & & \downarrow \\
\mathbf{Th}(-T\pi) & \longrightarrow & \mathbf{Th}(-\nu^*(\gamma)),
\end{array}$$

where $\underline{\mathbb{R}^2}$ denotes the trivial 2-dimensional vector bundle over $\overline{\mathrm{B}S\Gamma_2}$. From the diagrams 3.36 and 3.37, we obtain a homotopy commutative diagram

$$(3.38) \quad \begin{array}{ccc}
\mathrm{BDiff}_c^\delta(\mathbb{R}^2) & \xrightarrow{\alpha_{\mathbb{R}^2}} & Q_0 S^{-2}(\overline{\mathrm{B}S\Gamma_2}_+) \\
\downarrow & & \downarrow \\
\mathrm{BDiff}^\delta(\Sigma_g) & \xrightarrow{\alpha_{\Sigma_g}} & \Omega_0^\infty \mathrm{MT}\nu.
\end{array}$$

We are left to show that the map $\alpha_{\mathbb{R}^2}$ equals $\lambda \circ \beta$ up to homotopy. Recall that $\mathbb{R}^2 /\!\!/ \mathrm{Diff}_c^\delta(\mathbb{R}^2) \simeq \mathbb{R}^2 \times \mathrm{BDiff}_c^\delta(\mathbb{R}^2)$ is topologically trivial bundle and the flatness (transverse foliation on the total space) of the bundle gives rise to a map

$$f : \mathbb{R}^2 \times \mathrm{BDiff}_c^\delta(\mathbb{R}^2) \to \overline{\mathrm{B}S\Gamma_2},$$

that classifies the vertical tangent bundle $T\pi'$. Since the foliation on $\mathbb{R}^2 \times \mathrm{BDiff}_c^\delta(\mathbb{R}^2)$ is trivial outside of a compact set of the fiber, the map $f$ factors through the map $\mathbb{R}^2 \times \mathrm{BDiff}_c^\delta(\mathbb{R}^2) \to \Sigma^2(\mathrm{BDiff}_c^\delta(\mathbb{R}^2)_+)$. Therefore the spectrum map

$$\Sigma^\infty(\mathrm{BDiff}_c^\delta(\mathbb{R}^2)_+) \to \mathbf{Th}(-\underline{\mathbb{R}^2} \to \mathbb{R}^2 \times \mathrm{BDiff}_c^\delta(\mathbb{R}^2)) \to \mathbf{Th}(-\underline{\mathbb{R}^2} \to \overline{\mathrm{B}S\Gamma_2}),$$

whose adjoint is $\alpha_{\mathbb{R}^2}$ factors as

$$\Sigma^\infty(\mathrm{BDiff}_c^\delta(\mathbb{R}^2)_+) \xrightarrow{\mathrm{id}} \Sigma^{\infty-2}\Sigma^2(\mathrm{BDiff}_c^\delta(\mathbb{R}^2)_+) \to \mathbf{Th}(-\underline{\mathbb{R}^2} \to \overline{\mathrm{B}S\Gamma_2}).$$

Using $\Omega^\infty$-$\Sigma^\infty$ adjointness, the map $\alpha_{\mathbb{R}^2}$

$$\mathrm{BDiff}_c^\delta(\mathbb{R}^2) \to Q(\mathrm{BDiff}_c^\delta(\mathbb{R}^2)_+) \to QS^{-2}(\overline{\mathrm{B S\Gamma}_2}_+),$$

can be factored as

(3.39)
$$
\begin{array}{ccc}
\mathrm{BDiff}_c^\delta(\mathbb{R}^2) & \longrightarrow & Q(\mathrm{BDiff}_c^\delta(\mathbb{R}^2)_+) \\
{\scriptstyle \beta}\downarrow & & \downarrow \\
\Omega^2\overline{\mathrm{B S\Gamma}_2} & \xrightarrow{\ \lambda\ } & QS^{-2}(\overline{\mathrm{B S\Gamma}_2}_+).
\end{array}
$$

Hence, $\alpha_{\mathbb{R}^2} \simeq \lambda \circ \beta$. $\qquad\qquad\qquad\qquad\qquad\qquad\qquad\qquad\qquad\square$

Using the same idea, we will show that up to torsion, every codimension 2 foliation on a manifold of dimension 5 is foliated cobordant to a flat surface bundles of genus higher than 5.

Any flat $\Sigma_g$-bundle $[\Sigma_g \to E^{n+2} \xrightarrow{\pi} B^n]$ over an $n$-manifold $B^n$ gives a codimension 2 foliation on $E^{n+2}$. We let $e_n$ be the map that assigns to a flat surface the foliated cobordism class of the codimension 2 foliation on the total space. Hence, $e_n$ induces a well-defined map

$$e_n \colon \mathrm{MSO}_n(\mathrm{BDiff}^\delta(\Sigma_g)) \longrightarrow \mathcal{F}\Omega_{n+2,2}.$$

Let $E^5 \xrightarrow{\pi} B^3$ be a flat $\Sigma_g$-bundle over a 3 dimensional manifold $B^3$, then characteristic classes $h_1.c_2, h_1.c_1^2$ of the codimension 2 foliation on $E^5$ live in $H^5(E^5; \mathbb{R})$, hence

$$x = \langle \int_{\mathrm{fiber}} h_1.c_1^2, [B] \rangle, y = \langle \int_{\mathrm{fiber}} h_1.c_2, [B] \rangle$$

are real characteristic numbers associated to $E^5 \xrightarrow{\pi} B^3$. Consider the following diagram

$$
\begin{array}{ccc}
 & \mathcal{F}\Omega_{5,2} & \\
 & {\scriptstyle e_3}\nearrow \qquad \searrow {\scriptstyle (\int h_1 c_2, \int h_1 c_1^2)} & \\
\mathrm{MSO}_3(\coprod_g \mathrm{BDiff}^\delta(\Sigma_{g,1}, \partial)) & \xrightarrow{\ (x,y)\ } & \mathbb{R}^2.
\end{array}
$$

Since the classes $x, y$ are invariants of the foliated cobordism class, it is easy to see that the map induced by $(x, y)$ factors through $\mathcal{F}\Omega_{5,2}$. The surjectivity of the map induced by the integration of $h_1 c_2$ and $h_1 c_1^2$ is the statement of Rasmussen's theorem. In Theorem 3.27, we showed that the map induced by $(x, y)$ is rationally surjective. Now we prove that in fact $e_3$ is also rationally surjective.

**Theorem 3.40.** *The map*

$$e_3 \colon \mathrm{MSO}_3(\mathrm{BDiff}^\delta(\Sigma_g)) \longrightarrow \mathcal{F}\Omega_{5,2}$$

*is rationally surjective as $g \geq 5$, and is rationally isomorphism as $g \geq 6$.*

*Proof.* Recall that by a result of [Kos74, Theorem 1'], the map

$$\mathcal{F}\Omega_{k,2} \to \mathrm{MSO}_k(\mathrm{B S\Gamma}_2)$$

is rationally bijective for $k > 2$. Furhermore, by Theorem 2.6 and Theorem 1.14, we know that for $n \leq (2g-2)/3$, the following map

$$\mathrm{MSO}_n(\mathrm{BDiff}^\delta(\Sigma_g)) \longrightarrow \mathrm{MSO}_n(\Omega_0^\infty \mathrm{MT}\nu)$$

is bijective and for $n \leq 2g/3$, it is surjective. Consider the commutative diagram

(3.41)
$$\begin{array}{ccc} \mathrm{MSO}_3(\mathrm{BDiff}^\delta(\Sigma_g)) & \xrightarrow{\ e_3\ } & \mathcal{F}\Omega_{5,2} \\ \downarrow & & \downarrow \\ \mathrm{MSO}_3(\Omega_0^\infty \mathrm{MT}\nu) & \longrightarrow & \mathrm{MSO}_4(\mathrm{B}\Gamma_2). \end{array}$$

Therefore, the statement of the theorem is equivalent to proving that the bottom map in 3.41 is rationally an isomorphism. Note that the bottom map is given by composing the suspension map and the Thom isomorphism

$$\mathrm{MSO}_3(\Omega_0^\infty \mathrm{MT}\nu) \xrightarrow{\ \sigma\ } \mathrm{MSO}_3(\mathrm{MT}\nu) \xrightarrow{\text{Thom iso}} \mathrm{MSO}_5(\mathrm{B}\Gamma_2).$$

Recall that $\mathrm{MSO}_*(X) \otimes \mathbb{Q} \cong H_*(X;\mathbb{Q}) \otimes \mathrm{MSO}_*(\mathrm{pt})$ for any topological space $X$. Since both $\mathrm{B}\Gamma_2$ and $\Omega_0^\infty \mathrm{MT}\nu$ are simply connected one can easily see that

$$\mathrm{MSO}_3(\Omega_0^\infty \mathrm{MT}\nu) \otimes \mathbb{Q} \xrightarrow{\ \cong\ } H_3(\Omega_0^\infty \mathrm{MT}\nu;\mathbb{Q}),$$

$$\mathrm{MSO}_5(\mathrm{B}\Gamma_2) \otimes \mathbb{Q} \xrightarrow{\ \cong\ } H_5(\mathrm{B}\Gamma_2;\mathbb{Q}).$$

Therefore we need to show that the natural map

$$H_3(\Omega_0^\infty \mathrm{MT}\nu;\mathbb{Q}) \to H_3(\mathrm{MT}\nu;\mathbb{Q}),$$

is an isomorphism. We know from the proof of proposition 3.19 that the above is surjective. To prove that it is also injective consider the commutative diagram

$$\begin{array}{ccc} \pi_3(\Omega_0^\infty \mathrm{MT}\nu) \otimes \mathbb{Q} & \longrightarrow & \pi_3(\mathrm{MT}\nu) \otimes \mathbb{Q} \\ \downarrow & & \downarrow \\ H_3(\Omega_0^\infty \mathrm{MT}\nu;\mathbb{Q}) & \longrightarrow & H_3(\mathrm{MT}\nu;\mathbb{Q}). \end{array}$$

As we observed in proof of proposition 3.19, the composition of the top horizontal map and the right vertical map is an isomorphism. Given that the left vertical map is surjective by the rational Hurewicz theorem, the bottom horizontal map must be an isomorphism. □

Unlike $e_3$, the map

$$e_2 : \mathrm{MSO}_2(\mathrm{BDiff}^\delta(\Sigma_g)) \longrightarrow \mathcal{F}\Omega_{4,2}$$

is not rationally surjective. To every codimension 2 foliation $\mathcal{F}$ on a 4-manifold $M$, we can assign the difference of the Pontryagin classes $p_1(M) - p_1(\nu(\mathcal{F}))$, where $\nu(\mathcal{F})$ is the normal bundle of $\mathcal{F}$. It is easy to see that the number $\int_M p_1(M) - p_1(\nu(\mathcal{F}))$ is invariant of the foliated cobordism class of $\mathcal{F}$, hence it induces a map

$$\phi : \mathcal{F}\Omega_{4,2} \otimes \mathbb{Q} \to \mathbb{Q}.$$

Suppose we have a flat surface bundle $\Sigma_g \to M \xrightarrow{\ \pi\ } \Sigma_h$, then the normal bundle of the codimension 2 foliation on $M$ is the vertical tangent bundle $\mathrm{T}\pi$. It is easy to see that $p_1(M) = p_1(\mathrm{T}\pi) = p_1(\nu(\mathcal{F}))$ (see [Mor01, Proposition 4.11]). Hence $\phi$ vanishes on flat surface bundles. We prove below that vanishing of $\phi$ is essentially the only obstruction for a codimension 2 foliation $\mathcal{F}$ on a 4-manifold $M$ to be foliated cobordant to a flat surface bundle.

**Theorem 3.42.** *For $g \geq 4$, there is a short exact sequence*

(3.43)
$$0 \to \mathrm{MSO}_2(\mathrm{BDiff}^\delta(\Sigma_g)) \otimes \mathbb{Q} \xrightarrow{\ e_2\ } \mathcal{F}\Omega_{4,2} \otimes \mathbb{Q} \xrightarrow{\ \phi\ } \mathbb{Q} \to 0.$$

*Proof.* By [Kos74, Theorem 1'], we know that

$$\mathcal{F}\Omega_{4,2} \otimes \mathbb{Q} \xrightarrow{\cong} \mathrm{MSO}_4(\mathrm{B}\Gamma_2) \otimes \mathbb{Q}.$$

To show that $\phi$ is surjective, we just need to find a 4-manifold $M$ and a map $f : M \to \mathrm{B}\Gamma_2$ such that $\phi$ does not vanish on $[M, f] \in \mathrm{MSO}_4(\mathrm{B}\Gamma_2)$. Let $M$ be $\mathbb{C}P^2$ which is a 4 manifold whose signature is not zero, and $f$ be a nullhomotopic map. Since $f$ is trivial the normal bundle of the Haefliger structure $\mathcal{H}$ induced by $f$ is trivial. Thus, $p_1(\nu(\mathcal{H})) = 0$. Hence, we have

$$\phi([\mathbb{C}P^2, f]) = p_1(\mathbb{C}P^2) - p_1(\nu(\mathcal{H})) = 3 \neq 0.$$

To prove injectivity of $e_2$, note that we have

$$\mathrm{MSO}_2(\Omega_0^\infty \mathrm{MT}\nu) \otimes \mathbb{Q} \xrightarrow{\cong} H_2(\Omega_0^\infty \mathrm{MT}\nu; \mathbb{Q}) \xrightarrow{\cong} H_2(\mathrm{MT}\nu; \mathbb{Q}) \xrightarrow{\cong} H_4(\mathrm{B}\Gamma_2; \mathbb{Q}),$$

where the second isomorphism is given by the Hurewicz theorem because $\pi_1(\Omega_0^\infty \mathrm{MT}\nu) = 0$. Hence, using Theorem 2.6 and Theorem 1.14 we obtain for $g \geq 4$ the map

$$\mathrm{MSO}_2(\mathrm{BDiff}^\delta(\Sigma_g)) \otimes \mathbb{Q} \longrightarrow H_4(\mathrm{B}\Gamma_2; \mathbb{Q}),$$

is an isomorphism. On the other hand, by the Atiyah-Hirzebruch spectral sequence, we have a short exact sequence

$$0 \to \mathbb{Q} \to \mathrm{MSO}_4(\mathrm{B}\Gamma_2) \otimes \mathbb{Q} \to H_4(\mathrm{B}\Gamma_2; \mathbb{Q}) \to 0.$$

Therefore, we have a commutative diagram

$$
\begin{array}{ccc}
\mathrm{MSO}_2(\mathrm{BDiff}^\delta(\Sigma_g)) \otimes \mathbb{Q} & \xrightarrow{e_2} & \mathcal{F}\Omega_{4,2} \otimes \mathbb{Q} \\
\downarrow{\scriptstyle \cong} & & \downarrow{\scriptstyle \cong} \\
H_4(\mathrm{B}\Gamma_2; \mathbb{Q}) & \longleftarrow & \mathrm{MSO}_4(\mathrm{B}\Gamma_2) \otimes \mathbb{Q}.
\end{array}
$$

Hence, the map $e_2$ is injective with the cokernel $\mathbb{Q}$. Since $\phi$ vanishes on the image of $e_2$, the exactness in the middle term of 3.43 is also readily implied. $\square$

*Remark* 3.44. Jonathan Bowden pointed out to the author that in fact there is an example of a codimension 2 foliation (not just $S\Gamma_2$-structure) which is not in the image of $e_2$. To include his example, let $\mathcal{F}$ be the foliation by fibers of a surface bundle over a surface whose signature is nonzero. It is easy to see that $\phi(\mathcal{F}) \neq 0$, hence $\mathcal{F}$ is not in the image of $e_2$.

   *E-mail address*: sam@math.northwestern.edu

Department of Mathematics, Northwestern University, 2033 Sheridan Road,
Evanston, IL 60208